\documentclass{article}
\usepackage{my-settings, my-commands}
\usepackage{amsmath, amsfonts, amssymb}
\pdfoutput=1
\usepackage{graphicx}

\Testa

\newcounter{ContaLista}

\CadutaDelimitatori[1pt]

\def\K#1{\textit{#1}}
\def\T{\texttt}
\def\sm{\smallskip}
\def\NI{\noindent}
\def\Mp{\kern-3pt\longrightarrow\kern-2pt}
\def\Moid{\operatorname{id}}
\def\on{\operatorname}
\def\m{\medskip}
\def\Frac#1#2{\bgroup\displaystyle\frac{#1}{#2}\egroup}
\def\gk{\grf}
\def\odd{\on{odd}}
\def\Movg{\bigcup}
\def\epsilon{\varepsilon}
\def\Mkg{\gk}
\def\hme{^{-1}}
\def\Mkr{\tnd}
\newcommand\adjustqed[1][-14pt]{\belowdisplayskip=#1}
\def\Ls{\on{Ls}}
\def\Rs{\on{Rs}}

\def\U{\underline}

\def\gd{\sse}
\def\Mpi{\imp}
\def\Mf{\mathop{\bigcirc}\limits}

\newcommand\rahmen[2][10]{\spazio[-#1pt]\Allinea[16]{\boxed{\begin{array}{r}#2\end{array}}\notag}}

\def\Mke{\qdr}

\begin{document}
\bgroup
\Century[600]
\Titolo[]{\Century[1600]Dichotomic random number generators}{\Century[1100]Josef Eschgf\"{a}ller and Andrea Scarpante \\ \ \\ \Century[800] Universit\`a degli Studi di Ferrara \\ \ \\ \makeatletter esg@unife.it \& andrea.scarpante@student.unife.it \makeatother }
\egroup

\spazio[-10mm]

\Abstract{We introduce several classes of pseudorandom
sequences which represent a natural extension
of classical methods in random number
generation. The sequences are obtained from
constructions on labeled binary trees, generalizing
the well-known Stern-Brocot tree.}

\spazio[7mm]

\NI \textit{Keywords}: Dichotomic random number generator, pseudorandom sequence, binary tree, Stern-Brocot tree, Pari/GP.

\newpage

\Capitolo {Preliminaries}

\Situation {Let $X$ be a non-empty set.\sm

A \K{vector} is a finite (possibily void) sequence of elements of $X$.
In the combinatorics of words a vector is also called a (finite) \K{word} and the
set of all words is denoted by $X^*$. We shall use both terminologies.\sm

The length of a word $v$ is denoted by $|v|$.}

\Remark {For experiments, examples and graphical outputs we employed the computer algebra system Pari/GP, using
a collection of functions we prepared which is available on
\K{felix.unife.it/++/paritools}. The names of all functions in this collection begin
with \T{t\_}.}

\Definition {Let $a=(a_1,\ldots,a_m)$ and $b=(b_1,\ldots,b_{m+1})$ be two vectors
with $|b|=|a|+1$. Their \K{interleave} (or \K{shuffle}) $a\downarrow b$ is the vector\m

$(b_1,a_1,b_2,a_2,b_3,\ldots,b_m,a_m,b_{m+1})$\m

\NI One has
\vspace{-3mm}
\Allinea {(a\downarrow b)_{2j} &= a_j \quad\text{ for } j=1,\ldots,m\notag\\
(a\downarrow b)_{2j+1} &= b_{j+1}\quad\text{ for } j=0,\ldots,m\notag}}
\vspace{-3mm}
\Definition {The \K{natural binary tree} (NBT) is the infinite binary tree labeled
	by the elements of $\NN+1$ as in the figure:

\begin{figure}[h]
	\centering
	\includegraphics[width=0.7\linewidth]{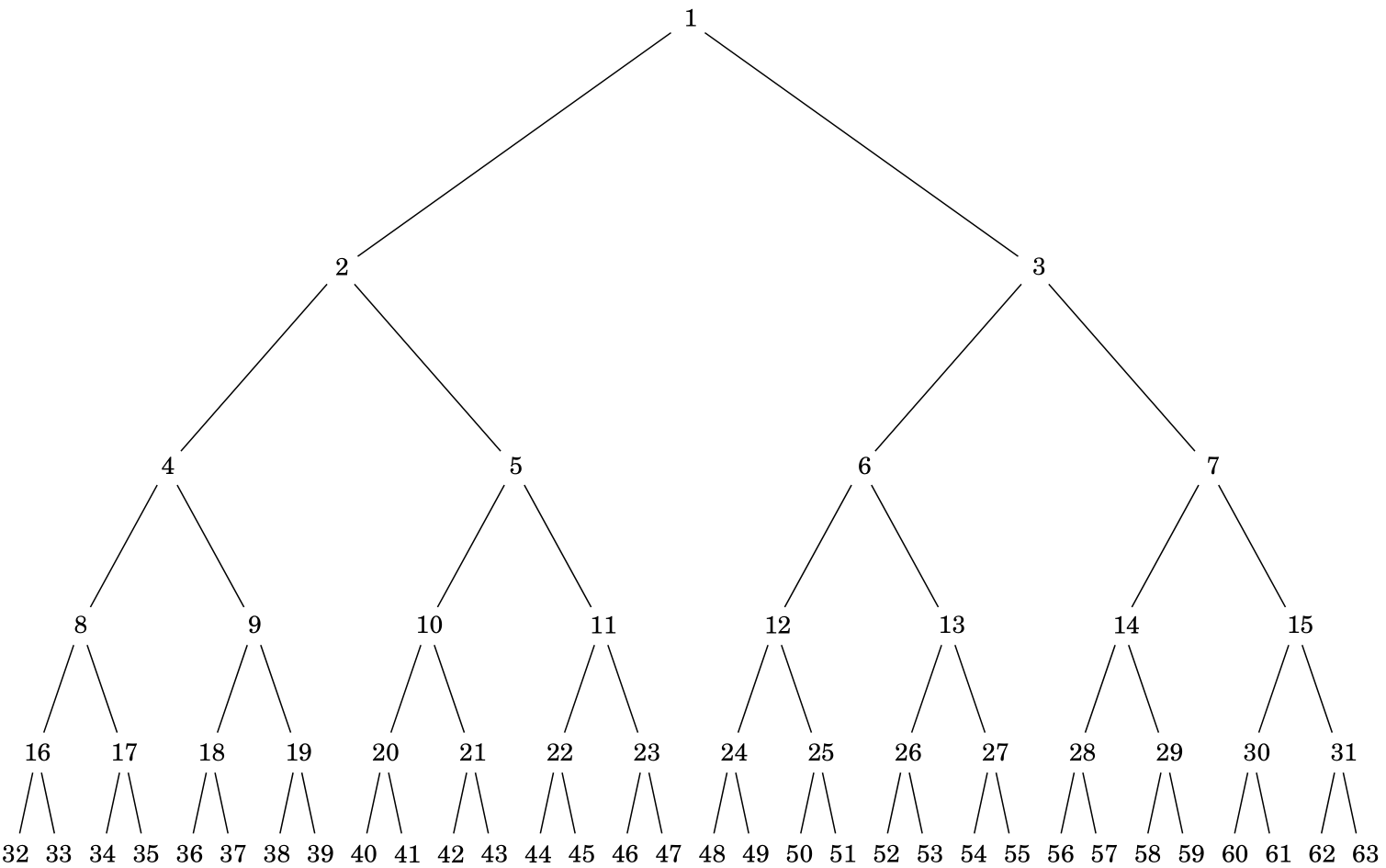}
\end{figure}

\NI The rows (row vectors)
of the tree are called its \K{levels}, the $k$-th level 
(beginning to count with $0$)
being denoted by $\LLL(*,k)$.}

\Remark {If $g:\NN+1\Mp X$ is a function, we obtain a labeled
tree $\LLL(g)$ whose levels are denoted
by $\LLL(g,k)$, as illustrated by the figure for \\ $g(n)=n^2$.

\begin{figure}[h]
	\centering
	\includegraphics[width=0.7\linewidth]{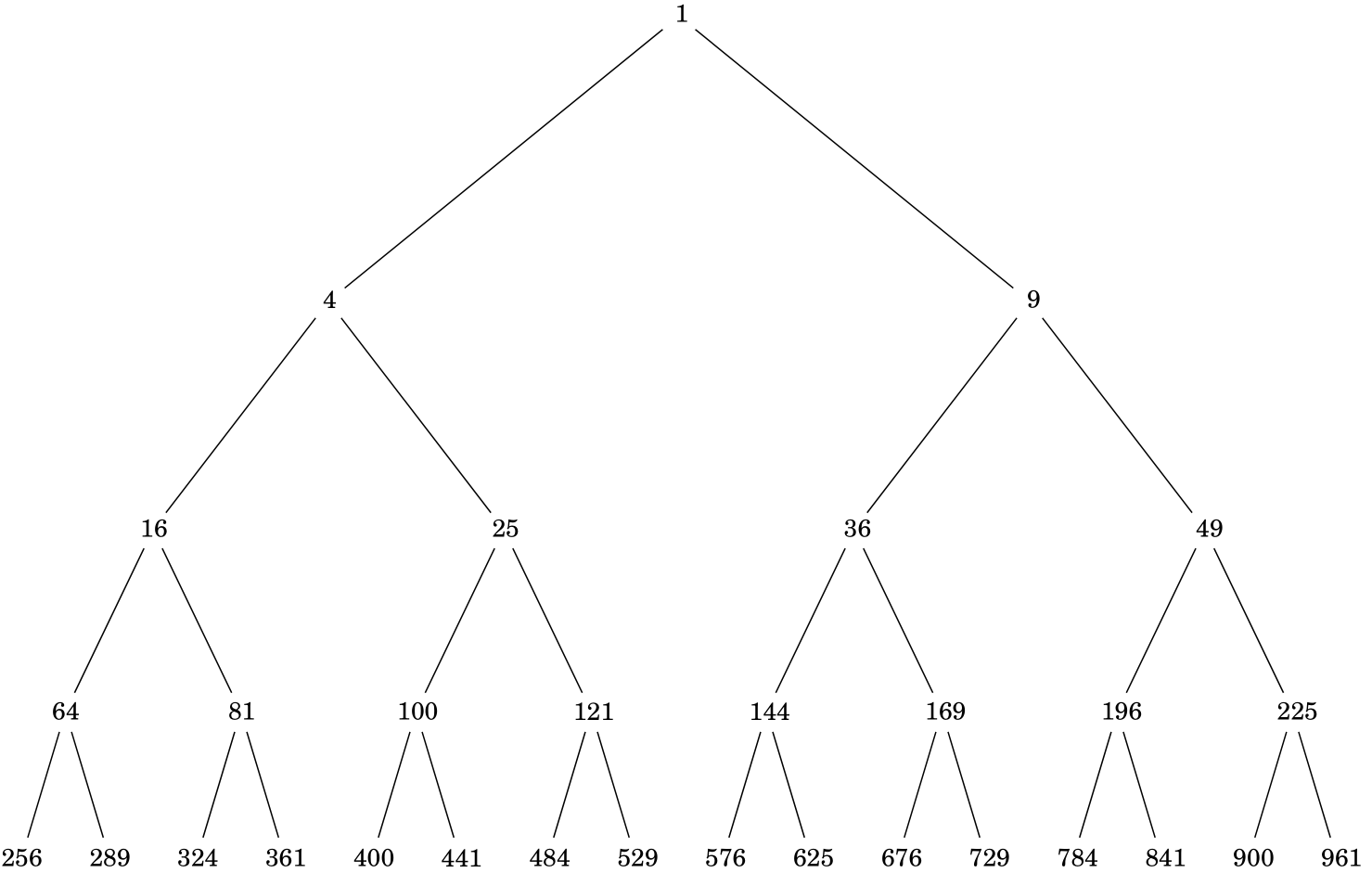}
\end{figure}

Hence $\LLL(*,k)=\LLL(\Moid,k)$, where $\Moid:\NN+1\Mp\NN+1$ is the identity function,
and the NBT can be written as $\LLL(*)$. 
More explicitly one has\sm

$\LLL(*,k)=(2^k,2^k+1,\ldots,2^{k+1}-1)$\sm

\NI and therefore\sm

$\LLL(g,k)=(g(2^k),g(2^k+1),\ldots,g(2^{k+1}-1))$\sm

\NI We count the elements in each row of the tree beginning with $1$ and denote
the $i$-th element of level $k$ by $\LLL(g,k,i)$. Hence\sm

$\LLL(g,k,i):= g(2^k+i-1)$}

\Definition {Every $n\in\NN+1$ belongs to a unique level $k$ and has therefore a
unique representation of the form\sm
	
$n=2^k+j$\sm
	
\NI with $k,j\in\NN$ and $0\le j<2^k$. In this case we write $n=2^k\oplus j$.\sm

We write also $L(n):=k$ for the level of $n$. Hence $j=n-2^{L(n)}$.\sm

In Pari/GP one obtains $L(n)$ as \T{\#binary(n)-1}.}

\Remark {We project now the NBT to the unit interval $[0,1]$ in such a way
that for $n=2^k\oplus j$ the abscissa $A(n)$ is given by\sm
	
$A(n)=\Frac{2j+1}{2^{k+1}}$\sm
	
\NI We obtain then a new labeled tree $\LLL(A)$, which is called the \K{dyadic tree}.
It contains every dyadic number $\Frac{2j+1}{2^{k+1}}$ with $k,j\in\NN$ and
$0\le j< 2^k$ exactly once.}

\begin{figure}[h]
	\centering
	\includegraphics[width=0.7\linewidth]{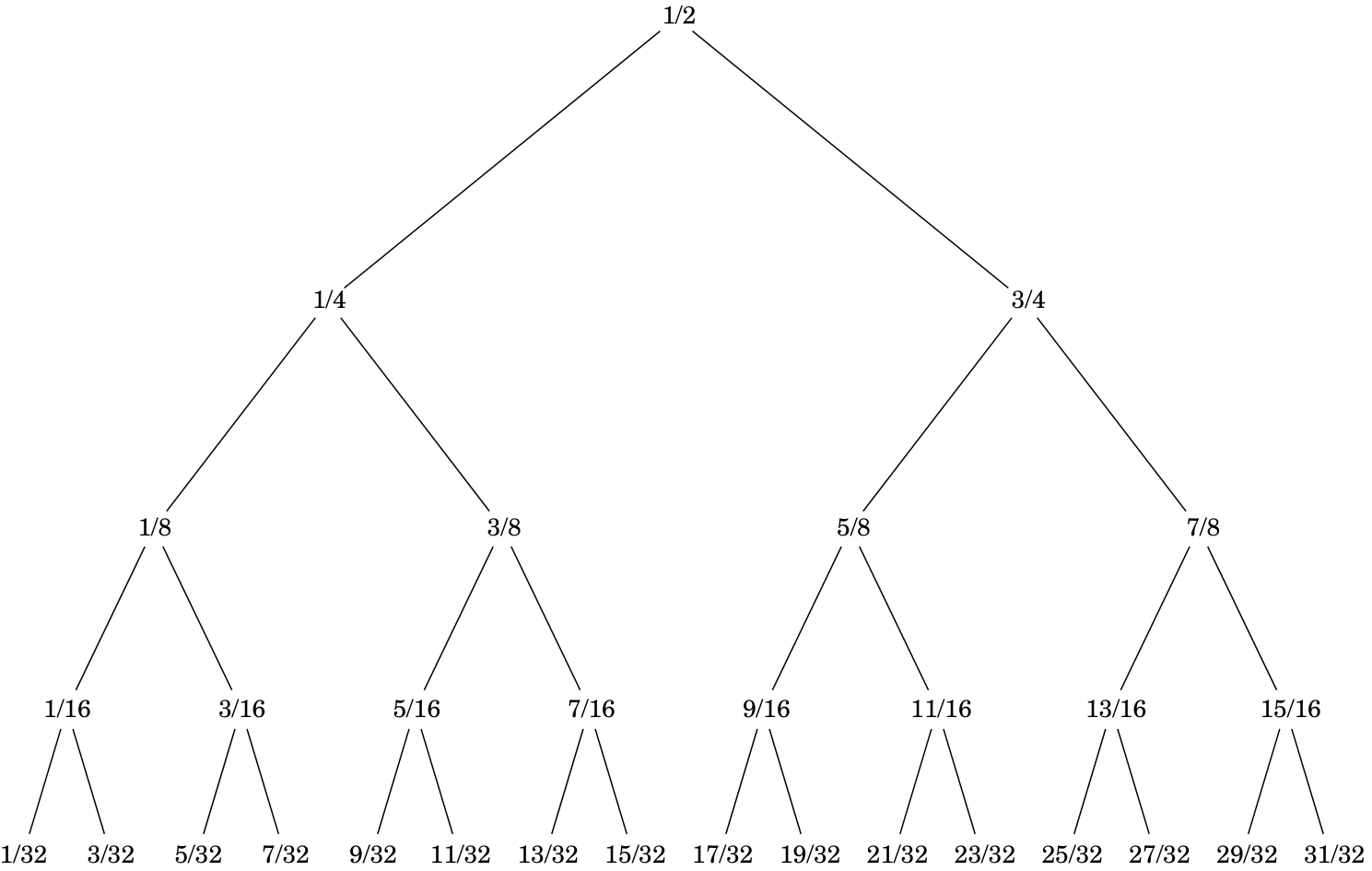}
\end{figure}

\Definition {Let $g:\NN+1\Mp X$ be a function and $S$ be a finite non-empty subset of
$\NN+1$. Assume that $S$ has exactly $m$ elements. Since the abscissa function $A$ of
Remark 1.7 is injective, we can write $S=\gk{s_1,\ldots,s_m}$ such that
$A(s_1)<A(s_2)<\ldots < A(s_m)$. See also Remark 1.17.\sm
	
The sequence $\EEE(g,S):= (g(s_1),\ldots,g(s_m))$ is then called the \K{binary
	evolution sequence} of $g$ on $S$.\sm
	
This is motivated by the following special case: For $k\in\NN$ let\\
$\NN(k):=\gk{n\in\NN+1 \mid n<2^{k+1}}$ be the full initial triangle up
to level $k$ of the NBT.
Then we can form the series of sequences
\Allinea[16]{\EEE(g,0) &:= \EEE(g,\NN(0)) = (g(1))\notag\\
	\EEE(g,1) &:= \EEE(g,\NN(1)) = (g(2),g(1),g(3))\notag\\
	\EEE(g,2) &:= \EEE(g,\NN(2)) = (g(4),g(2),g(5),g(1),g(6),g(3),g(7))\notag\\
	\ldots\notag}
	
\NI which is called the \K{binary evolution scheme} of $g$ and will
be denoted by $\EEE(g)$.\sm
	
We define $\EEE(g,-1)$ as the void sequence.\sm
	
Again we write $\EEE(*, \ldots)$ for $\EEE(\Moid,\ldots)$ and $\EEE(g,k,i)$ for
the $i$-th element\\of $\EEE(g,k)$. Hence\sm
	
$\EEE(g,k,i)=g(\EEE(*,k,i))$}

\Remark{In Def. 1.8 for every $k\in\NN+1$ one has\sm

$\EEE(g,k)=\EEE(g,k-1)\downarrow\LLL(g,k)$\sm

\NI From Definition 1.3 we have the recursion formulas\spazio
\Allinea[16]{\EEE(*,k,2j) &= \EEE(*,k-1,j) \quad\text{ for } j=1,\ldots,2^k-1\notag\\
	\EEE(*,k,2j+1) &= 2^k+j = \LLL(*,k,j+1) \quad\text { for } j=0,\ldots,2^k-1\notag}

\NI which in particular imply that\sm

$\EEE(*,k+\alpha,2^k) = 2^\alpha$ \quad for every $k,\alpha\in\NN$}

\Remark{The evolution scheme $\EEE(*)$ is interesting and well known:}
\vspace{-3mm}
\bgroup
\Century[850]
\begin{verbatim}
1                                                                                              
2   1   3                                                                                      
4   2   5   1   6   3   7                                                                      
8   4   9   2   10  5   11  1   12  6   13  3   14  7   15                                     
16  8   17  4   18  9   19  2   20  10  21  5   22  11  23  1   24  ... 
32  16  33  8   34  17  35  4   36  18  37  9   38  19  39  2   40  ... 
64  32  65  16  66  33  67  8   68  34  69  17  70  35  71  4   72  ... 
128 64  129 32  130 65  131 16  132 66  133 33  134 67  135 8   136 ... 
256 128 257 64  258 129 259 32  260 130 261 65  262 131 263 16  264 ... 
512 256 513 128 514 257 515 64  516 258 517 129 518 259 519 32  520 ... 
\end{verbatim}
\egroup

\NI or, if we want to respect the positions of the elements on the tree:
\bgroup
\Century[750]
\begin{verbatim}
                                        1
                  2                     1                     3
        4         2          5          1          6          3          7
   8    4    9    2    10    5    11    1    12    6    13    3    14    7    15
16 8 17 4 18 9 19 2 20 10 21 5 22 11 23 1 24 12 25 6 26 13 27 3 28 14 29 7 30 15 31
\end{verbatim}
\egroup
\NI Notice that $\EEE(*,k)$ is always a permutation of $\NN(k)$. This implies in particular
that $\EEE(*,k)$ has length $|\NN(k)|=2^{k+1}-1$.

Concatenating the vectors $\EEE(*,k)$ to an infinite sequence
$\EEE(*,0)\EEE(*,1)\cdots$, we obtain the sequence\sm

$(1,2,1,3,4,2,5,1,6,3,7,8,4,9,2,10,5,11,1,12,6,13,3,14,7,15,16,\ldots)$\sm

\NI which appears on OEIS as \K{A131987}. If one connects the same vectors by $0$, beginning
with $(0)$, one obtains the sequence\sm

$u=(0,0,1,0,2,1,3,0,4,2,5,1,6,3,7,0,8,4,9,2,10,5,11,1,12,6,13,3,\ldots)$\sm

\NI known as \K{A025480}. It is described by the simple recursion\sm

$u_{2n}=n, \quad u_{2n+1}=u_n$\sm

\NI beginning with $n=0$.\m

\Remark{ We observe first that the position in $\EEE(*,h)$
of a number $n$ which belongs to a level $\le h$ is given by
$A(n)\cdot 2^{h+1}$.\sm

If in the second output of Remark 1.10
we write only the new elements of each level,
we obtain a textual output of the NBT:}
\vspace{0mm}
\bgroup
\Century[750]
\begin{verbatim}
                                        1
                  2                                           3
        4                    5                     6                     7
   8         9         10         11         12         13         14         15
16   17   18   19   20    21   22    23   24    25   26    27   28    29   30    31
\end{verbatim}
\egroup

\Definition {We recall the following terminology from number theory:\sm
	
	Let $n\in\NN$. If $n>0$, then there exists a unique representation of the form
	$n=u\cdot2^m$ where $u$ is odd. We write $\on{odd}(n) := u$ and call it the \K{odd part} of
	$n$. Furthermore $|n|_2 := 2^{-m}$ is the $2$-\K{adic absolute value} of $n$.\sm
	
	We define $\on{odd}(0):=1$ and $|0|_2:=0$. Then:\sm
	
	(1) If $n>0$, then $\odd(n)$ is odd.\sm
	
	(2) $n$ is odd iff $\odd(n)=n$.\sm
	
	(3) $|n|_2 = 1$ iff $n$ is odd. In particular $|1|_2=1$.\sm
	
	(4) $\odd(n)=1$ iff $n=0$ or $n$ is a power of $2$.\sm
	
	(5) If $n>0$, then $n\cdot|n|_2=\odd(n)$.}

\Theorem {Let $k\in\NN$ and $0\le i<2^{k+1}$. Then}
	
	$\EEE(*,k,i) = 2^k|i|_2+\Frac{\odd(i)-1}{2}$
	
\Proof {Write $i=2^m\odd(i)$ and $\odd(i)=2j+1$. Then $m\le k$ and
	$0\le j<2^k$, so that from Remark 1.9 we obtain
	\Allinea[10]{\EEE(*,k,i) &= \EEE(*,k-m,\odd(i)) = \EEE(*,k-m,2j+1)\notag\\
		&= 2^{k-m}+j=2^{k-m}+\Frac{\odd(i)-1}{2}\notag}
	
	\NI Since $2^{-m}=|i|_2$, the theorem follows.}

\Corollary {Let $k,j\in\NN$ and $0\le j<2^k$. Then:}
	
	(1) $\EEE(*,k,2^k+j) = \Frac{2^{k+1}+2^k+j}{2}|j|_2-\Frac{1}{2}$.
	
	(2) If $j$ is odd, then $\EEE(*,k,2^k+j)=2^k+2^{k-1}+\Frac{j-1}{2}$.
	
\Proof {(1) The hypotheses on $j$ and $k$ imply that $|2^k+j|_2=|j|_2$.\\Since $2^k+j>0$, from Theorem 1.13 we have
	\Allinea {\EEE(*,k,2^k+j) &= 2^k|2^k+j|_2+\Frac{\odd(2^k+j)-1}{2}\notag\\
		&=2^k|2^k+j|_2+\Frac{(2^k+j)|2^k+j|_2-1}{2}\notag\\
		&= 2^k|j|_2+\Frac{(2^k+j)|j|_2-1}{2}=\Frac{2^{k+1}+2^k+j}{2}|j|_2-\Frac{1}{2}\notag}
	
	(2) This is a special case of (1) or also of Remark 1.9.}

\Proposition {Let $k,h\in\NN$ with $h\ge k$ and $n=2^k\oplus j\in\LLL(*,k)$. Then}\sm
	
	$n=\EEE(*,k,2j+1) = \EEE(*,h,(2j+1)\cdot 2^{h-k})$

\Proof {Immediate from Remark 1.9.}

\Remark{If we represent the NBT $\LLL(*)$ simply by its rows, we obtain the scheme}
\vspace{-3mm}
\bgroup
\Century[850]
\begin{verbatim}
1
2     3
4     5   6   7
8     9  10  11  12  13  14  15
16   17  18  19  20  21  22  23  24  25  26  27  28  29  30  31
32   33  34  35  36  37  38  39  40  41  42  43  44  45  46  47  48  49 ...
64   65  66  67  68  69  70  71  72  73  74  75  76  77  78  79  80  81 ...
128 129 130 131 132 133 134 135 136 137 138 139 140 141 142 143 144 145 ...
256 257 258 259 260 261 262 263 264 265 266 267 268 269 270 271 272 273 ...
512 513 514 515 516 517 518 519 520 521 522 523 524 525 526 527 528 529 ...
\end{verbatim}
\egroup

\NI The columns which appear in $\LLL(*)$ coincide with the columns which
appear in the scheme $\EEE(*)$ shown in Remark 1.10.

\Proof {This is immediate from Remark 1.9:\sm
	
	(1) Fix $i\in\NN+1$ and set $j:=i-1$. Then the $i$-th column in $\LLL(*)$ consists
	of the numbers $\LLL(*,k,i)$ with $k\in\NN$ such that $i\le 2^{k+1}$, i.e.
	$j<2^{k+1}$. By Remark 1.9 we have\sm
	
	$\LLL(*,k,i)=\LLL(*,k,j+1)=\EEE(*,k,2j+1)=\EEE(*,k,2i-1)$\sm
	
	(2) Fix again $i\in\NN+1$ and write $i=2^m(2j+1)$ with $m,j\in\NN$. As in the
	proof of Theorem 1.13 we have\sm
	
	$\EEE(*,k,i)=\EEE(*,k-m,2j+1)=\LLL(*,k-m,j+1)$}

\Remark [A very general method] {1. Let $(M,\prec)$ be
	totally ordered set and $g:M\Mp X$ be a
	mapping. Then each finite non-empty subset $S\subset M$ can be written in the
	form $S=\gk{s_1,\ldots,s_m}$ where $s_1\prec s_2 \prec \ldots \prec s_m$,
	giving rise to the vector $(g(s_1),\ldots,g(s_m))$. In some cases one could consider this
	vector as a \K{pseudorandom} sequence.\sm
	
	2. We shall apply this idea to the case $M=\NN+1$ and\sm
	
	$n\prec m \sse A(n)<A(m)$\sm
	
	\NI where $A$ is defined as in Remark 1.7. This order is known as \K{inorder} in
	computer science; cfr. Knuth \cita [p. 316-317] {12}. The sets $S$ will be often
	the sets $\NN(k)$ - the sequences generated are then the
	rows $\EEE(g,k)$ of the binary evolution scheme of $g$.\sm
	
	3. It could be interesting also to work with other subsets $S\subset \NN+1$.}

\newpage

\Capitolo {Generalized Stern-Brocot trees}	
	
\Situation {Let $X$ be a non-empty set. We use the standard notations
	from combinatorics of words:\vspace{-2mm}
	\Allinea[10] {X^* &:= \Movg_{n=0}^\infty X^n\notag\\
		X^+ &:= X^*\setminus\epsilon =  \Movg_{n=1}^\infty X^n\notag}
	
	\NI where $\epsilon$ is the empty word. Every $v\in X^*$ belongs to exactly one
	$X^n$ and we define then the length of $v$ as $|v|:= n$. In particular
	$|\epsilon|=0$.}

\Definition {Let $\overline{\NN} := \NN\cup\gk{1/2}$.
	
	We extend now the function $A$ of Remark 1.7 to a function
	$\overline{\NN}\Mp[0,1]$ by defining
	\Allinea[12]{ A(0) &:= 0\notag\\A(1/2) &:=1\notag}
	
	\NI The artificial elements $0$ and $1/2$ belong, by definition, to level $-1$.
	We put therefore $\LLL(*,-1) := (0,1/2)$
	
	Similarly we put, for any function $g: \overline{\NN}\Mp X$ and $k\in\NN$\sm
	
	$\overline{\EEE(g,k)} := g(0)\EEE(g,k)g(1/2)$\sm
	
	\NI and, as usual, $\overline{\EEE(*,k)} := \overline{\EEE(\Moid,k)}$.
	
	We shall not use the expressions $\overline{\EEE(g,k,i)}$, but define instead
	\Allinea[20] {\EEE(g,k,0) &:= g(0)\notag\\
		\EEE(g,k,2^{k+1}) &:= g(1/2)\notag}}

\Definition {Let $\DDD:=\Mkg{\Frac{a}{2^k} \mid a,k\in\NN \text{ with } 0<a<2^k}$ be
	the set of\\dyadic numbers and put\sm
	
	$\overline{\DDD} := \DDD\cup\gk{0,1} =
	\Mkg{\Frac{a}{2^k} \mid a,k\in\NN \text{ with } 0\le a\le 2^k}$\sm
	
	\NI The mapping $A$ from Remark 1.7 can then be considered as a mapping:\sm
	 
	$A:\overline{\NN}\Mp\overline{\DDD}$\sm
	
	with $A(0):=0$ and $A(1/2):=1$.\sm
	
	Notice that this mapping is bijective by construction.}

\Remark {Let $k\in\NN$ and $0\le i<2^{k+1}$. Then
	$A(\EEE(*,k,i))=\Frac{i}{2^{k+1}}$.}
\vspace{-3mm}
\Proof {Clear, since the projections of the elements of $\NN(k)$ are separated by intervals
	of length $\Frac{1}{2^{k+1}}$.
	
	Observe that the equation is true also for $i=0$, since $\EEE(*,k,0)=0$.}

\Proposition {Let $a,k\in\NN$ with $a\le 2^k$. Then}
	
	$A\hme \Mkr{\Frac{a}{2^k}} = 2^{k-1}|a|_2+\Frac{\odd(a)-1}{2}$

\Proof {(1) Consider first the case $0<a<2^k$. Then $\Frac{a}{2^k}
	\overset{2.4}{=} A(\EEE(*,k-1,a))$, hence\sm
	
	$A\hme\Mkr{\Frac{a}{2^k}}=\EEE(*,k-1,a)\overset{1.13}{=}
	2^{k-1}|a|_2+\Frac{\odd(a)-1}{2}$\sm
	
	(2) If $a=0$, then $2^{k-1}|a|_2+\Frac{\odd(a)-1}{2}=0=A\hme(0)$.\sm
	
	(3) If $a=2^k$, then $2^{k-1}|a|_2+\Frac{\odd(a)-1}{2}=\Frac{1}{2}+0
	=\Frac{1}{2}=A\hme(1)$.}

\Proposition {Let $k\in\NN$ and $1\le i<2^{k+1}$. If $i$ is odd, then
	\Allinea[20]{\EEE(*,k,i) &\ge2\cdot\EEE(*,k,i-1)+1\notag\\
		\EEE(*,k,i) &\ge2\cdot\EEE(*,k,i+1)\notag}}

\Proof {\adjustqed Since $i$ is odd, we have $|i|_2=1$ and $|i\pm1|_2\le\Frac{1}{2}$ and
	also $\odd(i\pm1)\le\Frac{i\pm1}{2}$. Writing
	for the moment $e_j:=\EEE(*,k,j)$ (for fixed $k$), from Theorem 1.13 now follow
	\Allinea{e_i &= 2^k|i|_2+\Frac{\odd(i)-1}{2} = \Frac{2^{k+1}}{2}\notag\\
		e_{i-1} &= 2^k|i-1|_2+\Frac{\odd(i-1)-1}{2}\le\Frac{2^k+\frac{i-1}{2}-1}{2}\notag\\
		&=\Frac{2^{k+1}+i-1}{4}-\Frac{1}{2}=\Frac{e_i-1}{2}\notag\\
		e_{i+1} &= 2^k|i+1|_2+\Frac{\odd(i+1)-1}{2}\le \Frac{2^k+\frac{i+1}{2}-1}{2}\notag\\
		&=\Frac{2^{k+1}+i-1}{4}=\Frac{e_i}{2}\notag}}

\vspace{12pt}

\Definition {For $k\in\NN$, the sequence $\EEE(*,k)$ contains,
	as noticed in Remark 1.9, all elements of $\LLL(*,k)$ in their natural order,
	interspersed with the elements of $\EEE(*,k-1)$, these belonging to levels $<k$, as shown
	here for level $k=3$, where we appended the two artificial elements on both
	extremities:\m
	
	$\begin{array}{*{17}{c}}
	0 & \mathbf{8} & 4 & \mathbf{9} & 2 & \mathbf{10} & 5 & \mathbf{11} & 1 &
	\mathbf{12} & 6 & \mathbf{13} & 3 & \mathbf{14} & 7 & \mathbf{15} & 1/2
	\end{array}$\m
	
	\NI The elements of $\LLL(*,3)$ are shown in boldface type. Similarly for every $k\in\NN$ each
	number $n\in\LLL(*,k)$ has a left and a right neighbor in $\overline{\EEE(*,k)}$, which
	belong to levels $<k$ and are called the \K{left support} $\Ls(n)$ and the
	\K{right support} $\Rs(n)$ of $n$ respectively.\sm
	
	It is also clear (by the very construction of $A$ in Remark 1.7) that
	\Allinea[12] {A(\Ls(n)) &= A(n)-\Frac{1}{2^{k+1}}\notag\\
		A(\Rs(n)) &= A(n)+\Frac{1}{2^{k+1}}\notag}
	
	\NI Notice finally that, since every $n\in\NN+1$ belongs to a unique level $k$, the left and
	the right support of $n$ are well defined for every such $n$.}

\Remark {(1) For $n\in\NN$ we have:\m
	
	\begin{tabular}{r@{\; $=$ \;}ll}
		$\Ls(2n)$ & $\Ls(n)$ & if  \; $n>0$\\
		$\Ls(2n+1)$ & $n$\\
		$\Rs(2n)$ & $n$ & if  \; $n>0$\\
		$\Rs(2n+1)$ & $\Rs(n)$ & if  \; $n>0$\\
		$\Rs(n)$ & $\Ls(n+1)$ & if  \; $n+1$ is not a power of $2$
	\end{tabular}\m
	
	(2) Moreover:\m
	
	\begin{tabular}{r@{\; $=$ \;}ll}
		$\Ls(2^k)$ & $0$ & for  \; $k\in\NN$\\
		$\Rs(2^k-1)$ & $1/2$ & for  \; $k\in\NN+1$
	\end{tabular}\m
	
	(3) In particular $\Ls(1)=0$ and $\Rs(1)=1/2$.}

\Proof {This is clear from the NBT.}
	
\Proposition {Let $n\in\NN+1$. Then $\Ls(n) = \Frac{\odd(n)-1}{2}$.}

\Proof {Write $n=2^m\odd(n)$ with $\odd(n)=2i+1$. By Remark 2.8 then\sm
	
	$\Ls(n)=\Ls(2i+1) = i = \Frac{\odd(n)-1}{2}$}

\Remark {Let $n\in\NN+1$.\m
	
	(1) If $n$ is even, then $\Rs(n)=\Frac{n}{2}> 2\Ls(n)$, hence $n>4\Ls(n)$.
	
	(2) If $n$ is odd $>1$, then $\Ls(n)=\Frac{n-1}{2}\ge2\Rs(n)$, hence $n>4\Rs(n)$.}
\vspace{-2mm}
\Proof {(1) From Remark 2.8 we know that $\Rs(n)=\Frac{n}{2}$. Now $n$ is even,
	therefore $\odd(n)\le\Frac{n}{2}$. Hence\sm
	
	$\Ls(n)\overset{2.9}{=}\Frac{\odd(n)-1}{2}\le\Frac{\frac{n}{2}-1}{2}
	=\Frac{n}{4}-\Frac{1}{2}$\sm
	
	\NI thus\sm
	
	$\Frac{n}{4}\ge\Ls(n)+\Frac{1}{2}>\Ls(n)$\sm
	
	(2) From Remark 2.8 we know that $\Ls(n)=\Frac{n-1}{2}$.
	
	Suppose first that $n+1$ is not a power of $2$. Then\sm
	
	$\Rs(n)=\Ls(n+1)=\Frac{\odd(n+1)-1}{2}\le\Frac{\frac{n+1}{2}-1}{2} = \Frac{n-1}{4}$\sm
	
	\NI hence\sm
	
	$\Frac{n}{4}\ge\Rs(n)+\Frac{1}{4}>\Rs(n)$\sm
	
	\NI Otherwise, if $n+1$ is a power of $2$, then $\Rs(n)=\Frac{1}{2}
	<\Frac{3}{4}\le \Frac{n}{4}$, since $n\ge3$ by hypothesis.}

\Definition {Let $f:X\times X\Mp X$ be a mapping and $a,b\in X$. Then we define
	a mapping $g:=f_{ab}:\overline{\NN}\Mp X$ in the following way:
	\Allinea[8] {g(0) &:= a\notag\\g(1/2) &:= b\notag\\
		g(n) &:= f(g(\Ls(n)),g(\Rs(n))) \text{ for } n\in\NN+1\notag}
	
	\NI Since for $n\in\NN+1$ the levels of $\Ls(n)$ and $\Rs(n)$ are both strictly
	smaller than the level of $n$, the mapping $f_{ab}$ is well defined.\sm
	
	Notice that always $g(1) = f(a,b)$.\sm
	
	Substituting each $n\in\NN+1$ in the NBT by $f_{ab}(n)$, we obtain the labeled
	binary tree $\LLL(f_{ab})$ which can be considered
	as a \K{generalized Stern-Brocot} tree, as we shall see (Proposition 2.14).}

\Remark {Let $g:\overline{\NN}\Mp X$ be a function and $k\in\NN$. Then\sm
	
	$\overline{\EEE(g,k)} = \LLL(g,k)\downarrow\overline{\EEE(g,k-1)}$}

\Proof {This follows from Remark 1.9, because appending one element on each side
	of the shorter sequence in Definition 1.3 corresponds to reversing the order of
	the two sequences around the $\downarrow$ symbol.}

\Remark {Let $k\in\NN$ and $n\in\LLL(*,k)$. Recall from Definition 2.7 that
	$\Ls(n)$ and $\Rs(n)$ are the left and right neighbors of $n$ in $\overline{\EEE(*,k)}$
	and thus are neighbors of each other in $\overline{\EEE(*,k-1)}$.\sm
	
	Consider now any function $g:\overline{\NN}\Mp X$. Then again $g(\Ls(n))$ and
	$g(\Rs(n))$ are neighbors of each other in  $\overline{\EEE(g,k-1)}$ and
	$g(n)$ is inserted between them in  $\overline{\EEE(g,k)}$.\sm
	
	If follows that, if now $f:X\times X\Mp X$, $a,b\in X$ and $g:=f_{ab}$,
	then $g(n)$ is the value of $f$ evaluated on the left and right neighbors of $g(n)$
	in $\overline{\EEE(g,k)}$ (taken in the position determined by $n$
	if it appears more than once), which both can be calculated on a lower level.\sm
	
	From Remark 2.12 we see that the sequence $\overline{\EEE(f_{ab},k)}$ is obtained from
	$x:=\overline{\EEE(f_{ab},k-1)}$ by inserting between $x_i$ and $x_{i+1}$ the
	value $f(x_i,x_{i+1})$.}

\Proposition {The NBT itself can be considered as a generalized Stern-Brocot tree.} 
	
	\textit{For this we define} $f:\overline{\NN}\times\overline{\NN}\Mp\overline{\NN}$ \textit{by}\sm 
	
	$f (x,y) := \Sistema{2y & \text{ if } x<y \\ 2x+1 & \text{ if } x>y \\ 0 & \text{ otherwise }}$\sm 
	
	\NI \textit{and choose} $a:=0$, $b:=1/2$.\textit{Then} $f_{ab} (n) = n$ \textit{for every} $n\in\overline{\NN}$.
	
\Proof {Let $g:=f_{ab}$. By definition $g(0)=0$, $g(1/2)=1/2$.
	
	Suppose $n\in\NN+1$.
	Then $n\in\LLL(*,k)$ for some $k\in\NN$. We use Remark 2.10 and show the proposition
	by induction on $k$.\sm
	
	\U{$k=0$:} Then $n=1$. But $g(1)=f(0,1/2) = 2\cdot\Frac{1}{2} = 1$.\sm
	
	\U{$k-1 \Mp k$:} If $n$ is even, then $\Rs(n) = \Frac{n}{2} > \Ls(n)$, hence\sm
	
	$g(n) = f (g(\Ls(n)),g(n/2)) \overset{IND}{=} f(\Ls(n),n/2)
	=2\Frac{n}{2} = n$\sm
	
	\NI If $n$ is odd, then $\Ls(n)=\Frac{n-1}{2}>\Rs(n)$, hence\sm
	
	$g(n) = f(g((n-1)/2),g(\Rs(n))) \overset{IND}{=} f((n-1)/2,\Rs(n))
	=n$}

\Definition {Let $f:X\times X\Mp X$ be a mapping, and $a,b\in X$.\sm
	
	Then we may construct a mapping $f_{ab}:\overline{\NN}\Mp X$ as in
	Definition 2.11.\sm
	
	The triple $(f,a,b)$ is called a
	\K{dichotomic generator} or simply a generator (of random sequences).}

\Remark {Let $f:X\times X\Mp X$ be mapping and $a,b\in X$.\\ For $k\in\NN$ then
the sequence $\EEE(f_{ab},k) = (x_1,\ldots,x_{2^{k+1}-1})$ can be calculated by the
general recursion formulas in Remark 1.9, but, as a consequence of Remark 2.13, also by
the following algorithm wich we describe in Pari/GP and which justifies the name
\K{dichotomic generator}:}
\vspace{-3mm}
\bgroup
\Century[950]
\begin{verbatim}
	dicho (f,n,i,j,a,b) = {my (m,x);
	if (n==i,return(a), n==j, return(b)); m=(i+j)\2;
	x=f(a,b); if (n==m,return(x));
	if (n<m, dicho(f,n,i,m,a,x), dicho(f,n,m,j,x,b))}
	
	\\ Example:
	f (x,y) = (3*x+5*y+2)%7
	p=2^101; q=2^94
	
	v=[dicho(f,n,0,p,0,1) | n<-[q..q+40]]
	t_to(v,60,,"")
	\\ 06562615052102001446426543236352445110603
\end{verbatim}
\egroup
\NI Notice that we may use this algorithm for calculating far away elements of the sequence
$\EEE(f_{ab},k)$, as we did in this example, where, for $f(x,y)=(3x+5y+2)\mod 7$,
$a=2$, $b=3$, $k=100$, the elements $x_n$ are calculated for
$n=2^{94}, 2^{94}+1,\ldots,2^{94}+40$. This calculation is done directly on these
indices without the need for calculating the preceding elements.

\Remark{Each finite sequence
	$x_0=a,\ldots,x_m=b$ of distinct elements
	can be obtained by the method of Remark 2.16:
	We define $f(x_0,x_m):=x_{[m/2]}$ and
	similarly $f(x_i,x_j):=x_{[(i+j)/2]}$, wherever
	these indices appear; all other values of $f$
	can be chosen arbitrarily.
	
	For example the sequence $(x_0,\ldots,x_{11})$
	can be obtained as a dichotomic sequence
	if we define:\spazio[-18pt]
	
\Allinea[14]{f(x_0,x_{11}) &:= x_5\notag\\
	f(x_0,x_5) &:= x_2\notag\\
	f(x_5,x_{11}) &:= x_8\notag\\
	f(x_0,x_2) &:= x_1\notag\\
	f(x_2,x_5) &:= x_3\notag\\
	f(x_5,x_8) &:= x_6\notag\\
	f(x_8,x_{11}) &:= x_9\notag\\
	f(x_6,x_8) &:= x_7\notag\\
	f(x_3,x_5) &:= x_4\notag\\
	f(x_9,x_{11}) &:= x_{10}\notag}
	}

\Remark {As far as we know, the idea of using Remark 2.13 for the generation
	of random sequences appears in Centrella \cita{9797} (written under
	the supervision of J. E.) and Kreindl \cita{24739}.}

\newpage

\Capitolo {Continuative Mappings}

\Remark {Let $g:\NN+1\Mp X$ be a mapping. \sm
	
	We shall then consider the sequences
	$\EEE(g,k)$ as (finite) random sequences, in the spirit of Remark 1.17. \sm
	
	For applications where
	unpredictability of the generated sequences is desired, as for example in cryptology,
	it may be a pleasing aspect of the method that the sequences $\EEE(g,k)$ for different $k$
	can be rather unrelated. For theoretical investigations, however, also the case that
	$\EEE(g,k+1)$ is always a continuation of $\EEE(g,k)$, i.e., that $\EEE(g,k)$ is always
	a prefix of $\EEE(g,k+1)$, will be interesting.\sm
	
	We shall now consider the question, when this happens, if $g$ is of the form $f_{ab}$
	as in Definition 2.11.}  

\Definition {Let $g:\NN+1\Mp X$ be a mapping. We define an infinite sequence
	$\EEE(g,\infty):\NN+1\Mp X$ by setting\sm
	
	$\EEE(g,\infty,n):= \EEE(g,\infty)(n) := \EEE(g,k,n)$\sm
	
	\NI if $n\in\LLL(*,k)$. This sequence consists of the values of $g$
	on the bold numbers in the following scheme (see Remark 1.10):\m
	
	\bgroup
	\Century[650]
		$\begin{array}{*{31}{r@{\hspace{1.2ex}}}}
		\mathbf 1 \\
		2 & \mathbf 1 & \mathbf 3 \\
		4 & 2 & 5 & \mathbf 1 & \mathbf 6 & \mathbf 3 & \mathbf 7 \\
		8 & 4 & 9 & 2 & 10 & 5 & 11
		& \mathbf 1 & \mathbf {12} & \mathbf 6 & \mathbf {13} &
		\mathbf 3 & \mathbf {14} & \mathbf 7 & \mathbf {15} \\
		16 & 8 & 17 & 4 & 18 & 9
		& 19 & 2 & 20 & 10 & 21 & 5 & 22 & 11 & 23
		& \mathbf 1 & \mathbf {24} & \mathbf {12} & \mathbf {25} & \mathbf 6
		& \mathbf {26} & \mathbf {13} & \mathbf {27} & \mathbf 3
		& \mathbf {28} & \mathbf {14} & \mathbf {29} & \mathbf 7 &
		\mathbf {30} & \mathbf {15} & \mathbf {31} \\
		\end{array}$\egroup\m
	
	\NI The bold numbers themselves represent the sequence $\EEE(*,\infty)$.\sm
	
	The sequence $\EEE(g,\infty)$, always defined, is of course interesting only 
	if\\ $\EEE(g,k+1)$ is
	a continuation of $\EEE(g,k)$ for every $k\in\NN$.
	
	In this case the mapping $g$ is called \K{continuative}.\sm
	
	If $g$ is defined on some set containing $\NN+1$ (usually on $\NN$ or on $\overline{\NN}$),
	this means, by convention, that the restriction $g_{|\NN+1}$ is continuative.}

\Remark {Since for $k,j\in\NN+1$ one has $2j+1\in\LLL(*,k)$ iff $j\in\LLL(*,k-1)$,
	the recursion formulas of Remark 1.9 become now
	\Allinea[16] {\EEE(g,\infty,1) &= g(1)\notag\\
		\EEE(g,\infty,2j) &= \EEE(g,\infty,j) \quad\text { for } j\in\NN+1\notag\\
		\EEE(g,\infty,2j+1) &= g(2^k+j) \quad\text { for } k\in\NN \text{ and }
		2^{k-1}\le j<2^k\notag}}

\Proposition {Let $g:\NN+1\Mp X$ be a mapping. Then the following statements
	are equivalent:\sm
	
	(1) $g$ is continuative.\sm
	
	(2) $g$ is constant on each column of $\EEE(*)$.\sm
	
	(3) $g$ is constant on each column of $\LLL(*)$.\sm
	
	(4) $\EEE(g,\infty,2j+1) = g(2^k+j)$ \quad for every $k,j\in\NN$ with $j<2^k$.\sm
	
	(5) $g(2^k+j)=g(2^m+j)$ for every $k,m,j\in\NN$ such that $j<2^k\le2^m$.\sm
	
	(6) $g(n)=g(n+2^{L(n)}\cdot(2^r-1))$ for every $n\in\NN+1, r\in\NN$.}\sm
	
	\NI Here $L(n)$ is the level of $n$ as in Definition 1.6. The rows and columns of $\EEE(*)$
	were represented in Remark 1.10, those of $\LLL(*)$ in Remark 1.16.\sm
	
	The columns of $\LLL(*)$ appear also as leftward diagonals in the tree-like representation
	(that is, in the NBT), as in the figure:
	
	\m\begin{minipage}{\linewidth}
		\makebox[\linewidth]{
		\includegraphics[width=0.7\linewidth]{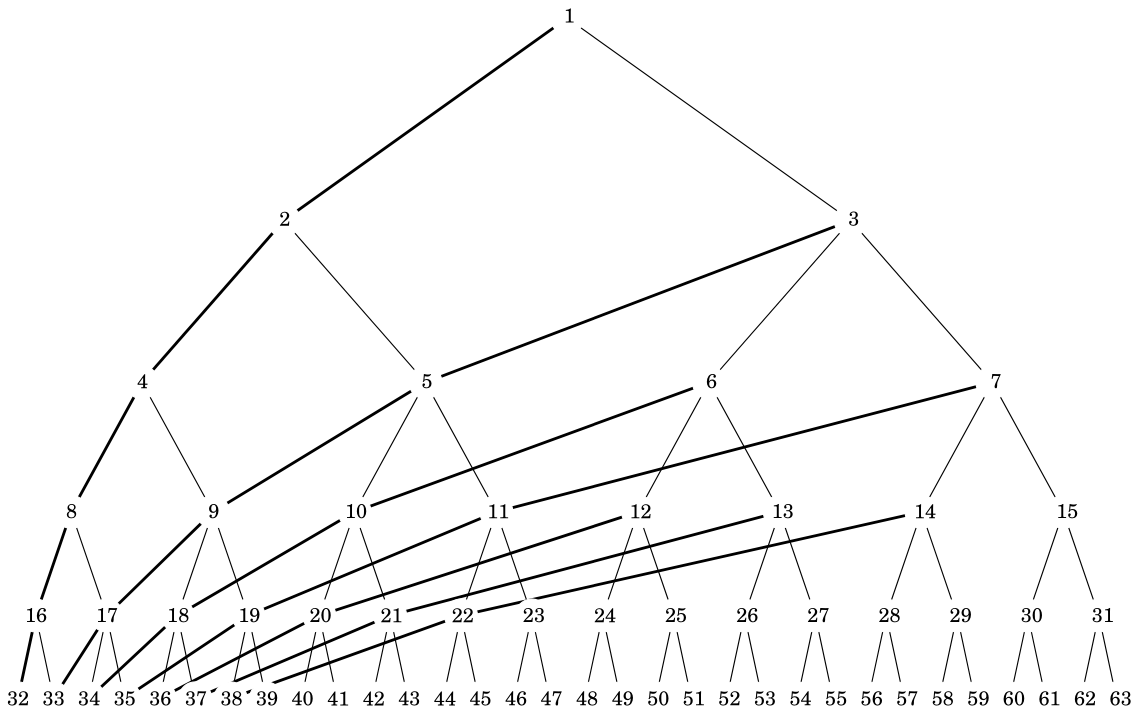}}
	\end{minipage}\m
	
	\Proof {(1) $\gd$ (2) $\gd$ (4) $\gd$ (5): By definition.\sm
		
		(2) $\gd$ (3): We observed in Remark 1.16 that $\LLL(*)$ and $\EEE(*)$ have the
		same columns - which in $\LLL(*)$ appear only once, in $\EEE(*)$ infinitely often.\sm
		
		(5) $\gd$ (6): Clear.} 
	
	\Lemma {Let $f:X\times X\Mp X$ be a mapping and $a,b\in X$. For every $k\in\NN$
		then\sm
		
		$\EEE(f_{ab},k+1) = \EEE(f_{a,f(a,b)},k)\cdot f(a,b)\cdot \EEE(f_{f(a,b),b},k)$\sm
		
		\NI where the dot denotes concatenation of words.}

	\Proof {Clear.}
	
	\Lemma {Let $f:X\times X\Mp X$ be a mapping and $b,c\in X$.
		Assume that $f(x,b)=f(x,c)$ for every $x\in X$.
		
		Then $f_{ab}=f_{ac}$ for every $a\in X$.}

	\Proof {We show by induction on $k\in\NN$ that $\EEE(f_{ab},k)
		=\EEE(f_{ac},k)$ for every $a\in X$ and every $k\in\NN$.\sm
		
		\U{$k=0$:} Applying the hypothesis to $x=a$ we have $f(a,b)=f(a,c)$, hence\sm
		
		$\EEE(f_{ab},0) = (f(a,b)) = (f(a,c)) = \EEE(f_{ac},0)$\sm
		
		\U{$k\Mp k+1$:} One has
		\Allinea[10] {\EEE(f_{ab},k+1) &\overset{3.5}{=}
			\EEE(f_{a,f(a,b)},k)\cdot f(a,b)\cdot\EEE(f_{f(a,b),b},k)\notag\\
			&= \EEE(f_{a,f(a,c)},k)\cdot f(a,c)\cdot\EEE(f_{f(a,c),b},k)\notag\\
			&\overset {IND}{=} \EEE(f_{a,f(a,c)},k)\cdot f(a,c)\cdot\EEE(f_{f(a,c),c},k)
			= \EEE(f_{ac},k+1)\notag}
		
		\NI where we used again that $f(a,b)=f(ac)$, applying in $\overset{IND}{=}$ the
		induction hypothesis on $f(a,c)$ instead of $a$.} 
	
\Proposition {Let $f:X\times X\Mp X$ be a mapping and $a,b\in X$.
	Assume that $f(x,f(a,b))=f(x,b)$ for every $x\in X$.
	
	Then $f_{ab}$ is continuative.}

\Proof {For every $k\in\NN$ we have\sm
	
	$\EEE(f_{ab},k+1) \overset{3.5}{=} \EEE(f_{a,f(a,b)},k)
	\cdot f(a,b)\cdot \EEE(f_{f(a,b),b},k)$ \mybox{(*)}\sm
	
	\NI The hypothesis $f(x,b)=f(x,f(a,b))$ for every $x\in X$ implies by Lemma 3.6 that
	$f_{a,f(a,b)}=f_{ab}$, hence (*) implies that 
	$\EEE(f_{ab},k)=\EEE(f_{a,f(a,b)},k)$ is a
	prefix of $\EEE(f_{ab},k+1)$.}

\Corollary {Let $f:X\times X\Mp X$ be a mapping and $a,b\in X$.
	
	If $f(a,b)=b$, then $f_{ab}$ is continuative.}

\Remark {Let $g:\overline{\NN}\Mp X$ be a mapping and set $a:=g(0)$, $b:=g(1/2)$.
	Consider the sequences $\overline{\EEE(g,k)}$:\m
	
	$\begin{array}{*9c}
	a & g(1) & b\\
	a & g(2) & g(1) & g(3) & b\\
	a & g(4) & g(2) & g(5) & g(1) & g(6) & g(3) & g(7) & b\\
	\ldots
	\end{array}$\sm
	
	\NI Then, for any fixed $k\in\NN$, $\overline{\EEE(g,k+1)}$ is a continuation
	of $\overline{\EEE(g,k)}$ iff $\EEE(g,k+1)$ is a continuation of $\EEE(g,k)$ and,
	in addition, $g(1)=b$.} 

\Corollary {Let $f:X\times X\Mp X$ be a mapping and $a,b\in X$. The following
	statements are equivalent:\sm
	
	(1) $\overline{\EEE(f_{ab},k+1)}$ is a continuation of
	$\overline{\EEE(f_{ab},k)}$ for every $k\in\NN$.\sm
	
	(2) $\EEE(f_{ab},k+1)$ is a continuation of $\EEE(f_{ab},k)$ for every $k\in\NN$ and,
	in addition, $f(a,b)=b$.\sm
	
	(3) $f(a,b)=b$.}

\Proof {(1) $\gd$ (2): Remark 3.9.\sm
	
	(2) $\Mpi$ (3): Clear.\sm
	
	(3) $\Mpi$ (2): Corollary 3.8.\sm
	
	\NI We found this result first in Kreindl \cita{24739}.} 	

\newpage

\Capitolo{One-sided generators}

\Situation {Let $X$ be a non-empty set.}

\Definition {If $P$ is a property defined for mappings, we say that the generator $(f,a,b)$ has
	property $P$ if the mapping $f_{ab}$ has property $P$. Thus for example
	the generator $(f,a,b)$
	is called continuative, if the mapping $f_{ab}$ is continuative.}

\Definition {A dichotomic generator $(f,a,b)$ is called \K{one-sided}, if
	$f(x,y)$ depends only on $x$. In this case there exists a function $\phi:X\Mp X$
	such that $f(x,y)=\phi(x)$ for every $x,y\in X$.}

\Remark {Every one-sided generator is continuative.}

\Proof {Let $(f,a,b)$ be a one-sided generator. 
	
	For every $x\in X$ then
	$f(x,f(a,b)) = f(x,b)$, since $f$ does not depend on the second argument. Hence $f_{ab}$ is continuative by Proposition 3.6.} 

\Definition {Let $\phi:X\Mp X$ be a mapping and $a\in X$. We define a mapping
	$g:\NN\Mp X$ in the following way:\vspace{-8pt}
	\Allinea[16] {g(0) &:= a\notag\\g(n) &:= \phi(g(\Ls(n))) \quad\text{ for } n\in\NN+1\notag}
	
	\NI and write also $\phi_a:=g$. Since for $n>0$ always $\Ls(n)<n$, the mapping
	is well defined.
	
	On its domain of definition $\phi_a$ coincides obviously with $f_{ab}$, if we define
	$f(x,y):=\phi(x)$ and choose $b\in X$ arbitrarily. Therefore we shall also call
	the couple $(\phi,a)$ or, for short, the mapping $\phi_a$ itself, a one-sided generator.} 

\Proposition {Let $\phi:X\Mp X$ be a mapping and $a\in X$. Then:\spazio
	\Allinea[16] {\phi_a(2j) &= \phi_a(j)\notag\\
		\phi_a(2j+1) &= \phi(\phi_a(j))\notag}
	
	\NI for every $j\in\NN$. In particular $\phi_a(1) = \phi(a)$.}

\Proof {(1) This statement is trivial for $j=0$. Assume $j>0$. Then\sm
	
	$\phi_a(2j)=\phi(\phi_a(\Ls(2j)))\overset{2.8}{=}\phi(\phi_a(\Ls(j)))=\phi_a(j)$.\m
	
	(2) $\phi_a(2j+1)=\phi(\phi_a(\Ls(2j+1)))\overset{2.8}{=}\phi(\phi_a(j))$.} 

\Theorem {Let $\phi:X\Mp X$ be a mapping and $a\in X$. Then\sm
	
	$a\EEE(\phi_a,\infty)=\phi_a$\sm
	
	\NI or, equivalently,\sm
	
	$\EEE(\phi_a,\infty,n) = \phi_a(n)$\sm
	
	\NI for every $n\in\NN+1$.}

\Proof {\adjustqed Let $u:=a\EEE(\phi_a,\infty)$, hence $u_0=a$ and
	$u_n=\EEE(\phi_a,\infty,n)$ for $n\in\NN+1$.\sm
	
	(1) We show that $u$ satisfies the same recursion rules as $\phi_a$, i.e., that\spazio
	\Allinea[-8] {u_1 &= \phi(a)\notag\\u_{2j} &= u_j\notag\\u_{2j+1} &= \phi(u_j)\notag}
	
	\NI for every $j\in\NN+1$. This clearly implies $u=\phi_a$.\sm
	
	(2) Since by Remark 4.4 $\phi_a$ is continuative, from Proposition 3.4 we have\spazio
	\Allinea[-8] {u_1 &= \phi_a(1) = \phi(a)\notag\\
		u_{2j} &= u_j\ \tn{ for } j\in\NN+1\notag\\
		u_{2j+1} &= \phi_a(2^k+j)\ \tn{ for every }k,j\in\NN\tn{ with }j<2^k\notag}
	
	(3) We show by induction on $k\in\NN$ the following statement:\sm
	
	If  $0\le j<2^k$, then $u_{2j+1}=\phi(u_j)$.\sm
	
	\U{$k=0$:} In this case $j=0$ and we have to show that $u_1=\phi(u_0) =\phi(a)$,
	and this is true.\sm
	
	\U{$k-1\Mp k$:} Assume $0\le j<2^k$.\sm
	
	Suppose first that $j$ is odd. Since now $k>0$, also $2^k+j$ is odd, thus\spazio
	\Allinea[-8] {u_{2j+1} &= \phi(\phi_a(\Ls(2^k+j)))\overset{2.8}{=}
		\phi\Mkr{\phi_a\Mkr{\Frac{2^k+j-1}{2}}}\notag\\
		&= \phi\Mkr{\phi_a\Mkr{2^{k-1}+\Frac{j-1}{2}}}\overset{(2)}{=}
		\phi\tnd{u_{2\frac{j-1}{2}+1}} = \phi(u_j)\notag}
	
	\NI since $0\le\Frac{j-1}{2}<2^{k-1}$.\sm
	
	Suppose now that $j$ is even. For $j=0$ we have $u_1=\phi(u_0)=\phi(a)$ as before.
	Otherwise write $j=2^mr$ with $r$ odd. Then $0<m<k$ and\spazio
	\Allinea {u_{2j+1} \overset{(2)}{=} \phi_a(2^k+j) = \phi_a(2^m(2^{k-m}+r))
		\overset{4.6}{=}\phi_a(2^{k-m}+r)\notag\\
		\overset{(2)}{=}u_{2r+1}\overset{IND}{=} \phi(u_r)=\phi\tnd{u_{\frac{j}{2^m}}}
		\overset{(2)}{=}\phi(u_j)\notag}}\sm 	
	
\Remark {The conclusion in Theorem 4.7 is not more true for general continuative
	dichotomic generators, as the example $(f,1,6)$ with\
	$f(x,y)=(3x+2y+7) \mod 8$ shows:}

\spazio[-3mm]	
\bgroup
\Century[950]
\begin{verbatim}
g=f_{1,6}     : 6 6 5 6 5 3 2 6 5 3 2 7 2 2 1 6 5 3 2 7 2 2 1 7 2 4 7 2 ...
E(g,infinite) : 6 6 5 6 3 5 2 6 7 3 2 5 2 2 1 6 7 7 2 3 4 2 7 5 2 2 1 2 ...
\end{verbatim}
\egroup\spazio[-8pt]

\Remark {Since in the proof of
	Theorem 4.7 $u$ is uniquely determined by the recursion rules (*),
	for a sequence $u\in X^\NN$ with $a:=u_0$ and a mapping $\phi:X\Mp X$
	the following statements are equivalent:\sm
	
	(1) $u=\phi_a$.\sm
	
	(2) $u=a\EEE(\phi_a,\infty)$.\sm
	
	(3) For every $j\in\NN$ we have $u_{2j}=u_j$ and $u_{2j+1}=\phi(u_j)$.\sm
	
	\NI The infinite sequences which obey a recursion rule of type (3) are therefore exactly the
	sequences obtained by a one-sided generator as in Theorem 4.7.} 

\EDS{Example} {Let $u$ be the Thue-Morse sequence $u\in\gk{0,1}^\NN$ defined by
	$u_0:= 0, u_{2j}=u_j, u_{2j+1}=1-u_j$.\sm
	
	By Theorem 4.7 it can be obtained as $u=0\EEE(\phi_0,\infty)$, where
	$\phi(x):=1-x$.}

\EDS{Example} {Consider the function $\beta:=\Mf_nn+1:\NN\Mp\NN$ and define\\
	$h:=\beta_0:\NN\Mp\NN$ in accordance with Definition 4.5 by
	\Allinea[12] {h(0) &:= 0\notag\\h(n) &:= h(\Ls(n))+1 \quad \text{ for } n\in\NN+1\notag}
	
	\NI Then by Theorem 4.7\sm
	
	$h=0\EEE(h,\infty) = 0\ 1\ 1\ 2\ 1\ 2\ 2\ 3\ 1\ 2\ 2\ 3\ 2\ 3\ 3\
	4\ 1\ 2\ 2\ 3\ 2\ 3\ 3\ 4\ 2\ 3\ \ldots$\sm
	
	\NI One can also show that $h(n)$ is equal to the Hamming weight of $n$, i.e. to the
	number of ones in the binary representation of $n$. This sequence is well known and listed
	as \K{A000120} in the OEIS.} 

\Proposition {Let $\phi:X\Mp X$ be a mapping and $a\in X$.
	
	Define $h$ as in Example 4.11. Then\sm
	
	$\phi_a(n) = \phi^{h(n)}(a)$\sm
	
	\NI for every $n\in\NN$.}

\Proof {\adjustqed We show the proposition by induction on $n$.\sm
	
	\U{$n=0$:} $\phi^{h(0)}(a)=\phi^0(a)=a=\phi_a(0)$.\sm
	
	\U{$n-1 \Mp n$:} Now $n>0$ and we may write $n=2^mr$ with $r$ odd.\sm
	
	By Proposition 2.9 $\Ls(n)=\Frac{r-1}{2}$, hence $h(n)=h\Mkr{\Frac{r-1}{2}}+1$. Further
	\Allinea {\phi_a(n) &= \phi_a(2^mr) \overset{4.6}{=}\phi_a(r)
		=\phi\Mkr{\phi_a\Mkr{\Frac{r-1}{2}}}\notag\\
		&\overset{IND}{=}\phi\Mkr{\phi^{h\tnd{\frac{r-1}{2}}}(a)}
		=\phi^{h\tnd{\frac{r-1}{2}}+1}(a) = \phi^{h(n)}(a)\notag}} 

\newpage

\Capitolo{Examples}

\Remark {In the following chapter we present examples of dichotomic generators.\sm
	
	For every generator are first indicated the function $f:X\times X\Mp X$
	(with $X$ usually tacitly understood) and the initial values $a,b\in X$.\sm
	
	Then follows the beginning of the binary evolution scheme (Definition 1.8) of
	the function $f_{ab}$, from which the last row is
	selected. This vector of values is represented graphically in a bar diagram;
	by a similar bar diagram we represent also the absolute values of
	the discrete Fourier transform of the vector, with the origin centered.\sm 
	
	Using the values $x_i-\mu$ as increments, where $\mu$ is the mean of
	the vector, we obtain a random walk which is given too.\sm 
	
	On the left then we present a usually longer vector of the same
	level of the evolution scheme by points in the plane, which are
	calculated in the following manner: As for the discrete Kolmogorov-Smirnov
	test (cf. Centrella \cita{9797}) first the vector  is decomposed in ordered non-overlapping blocks of length $10$. 
	Then the Ruffini-Horner method for powers of $2$
	is applied to each block giving us a vector of real numbers:\m
	
	$u:=(u_1,...,u_r)$\m
	 
	where $r$ is the number of blocks. 
	
	\NI Finally each entry $u_i$ of $u$ is divided 
	by $2^{10}$, which gives the vector\m
	
	$v:=(v_1,...,v_r)$ with $v_i:=\Frac{u_i}{2^{10}}$.\m
	
	\NI Now from the pairs $(v_{2k},v_{2k+1})$ we obtain a $2$-dimensional 
	representation of the sequence.
	
	%\sm Each time the numerical results of the battery of tests are given,
	%using the shortcuts defined in Remark 5.14.
	} 
	
\Remark{Each time the numerical results of a battery of tests are given using the following shortcuts:\bigskip
	
			\begin{tabular}{l @{\ \Puntini\ } l}
				\texttt{runs}  & \begin{minipage}{3.47in} \textit{run test}  (cf. Bassham a.o. \cita{25668}, Maurer \cita{25830} and Fisz \cita{2031})\end{minipage}\m \\ 
				\texttt{freq}  & \begin{minipage}{3.5in} \textit{frequency test}  (cf. Bassham a.o. \cita{25668})\end{minipage}\m \\ 
				\texttt{cusum}  & \begin{minipage}{3.5in} \textit{cumulative sum test}  (cf. Bassham a.o. \cita{25668})\end{minipage}\m \\ 
				\texttt{blocks}  & \begin{minipage}{3.5in} \textit{blocks test}  (cf. Fisz \cita{2031})\end{minipage}\m \\ 
				\texttt{autocorr} & \begin{minipage}{3.5in} \textit{auto-correlation test}  (cf. Bassham a.o. \cita{25668})\end{minipage}\m \\ 
				\texttt{longrun} & \begin{minipage}{3.5in} \textit{longrun test}  (cf. Guibas \& Odlyzko \cita[p. 252-253]{25885})\end{minipage}\m \\ 
				\texttt{2bits} & \begin{minipage}{3.5in} \textit{$2$-bit test}  (cf. Fisz \cita[page 399]{2031} and Bassham a.o. \cita{25668})\end{minipage}\m \\
				\texttt{ks\_discrete} &  \begin{minipage}{3.1in}  \textit{discrete Kolmogorov-Smirnov test}   (cf. Kuipers \& Niederreiter \cita[p. 90-92]{19061} and Fisz \cita{2031})\end{minipage} 
			\end{tabular}
			\newpage
			\begin{tabular}{l @{\ \Puntini\ } l}
				\texttt{DTF} & \begin{minipage}{3.5in} \textit{discrete Fourier transform test}   (cf. Bassham a.o. \cita{25668})\end{minipage}\m \\ 
				\texttt{Maurer} & \begin{minipage}{3.7in}
					\textit{Maurer's universal test}   (cf. Maurer \cita{25830},  Coron \& Naccache \cita{26401}, Do\u{g}anaksoy \& Tezcan \cita{26390} and Bassham a.o. \cita{25668})
				\end{minipage}
			\end{tabular}\bigskip		
	
\NI For every example we simply project the generate sequence onto $\ZZ/2\ZZ$ and we apply the above, most commonly used, bit tests, as described in the cited references. }

\newpage

\EDS{Example}{\spazio[-10pt]
\Allinea[16]{\boxed{\begin{array}{r}
	f(x,y) = (x+y+1) \mod 7  \\[2ex]   a=3, b=5
	\end{array}}\notag}   }
\bgroup
\Century[950]
\begin{verbatim}
2
621
3622410
033622520461206
4043033622521512305446114230065
145054134043033622521512410501426340653424461131640263401006554
5164356065346153145054134043033622521512410501420461206560216402263 ...
2501164413255600065523144611052351643560653461531450541340430336225 ...
6215602131164424615362154556001010065545126351642446113120651263250 ...
3622410556003241535131164424020446110523362241053435455600102120212 ...
0336225204612065455600104362046105232501535131164424020450323054244 ...
4043033622521512305446114230065534354556001021205413362230544611206 ...
\end{verbatim}
\egroup
\begin{figure}[h]
\includegraphics[width=1\linewidth]{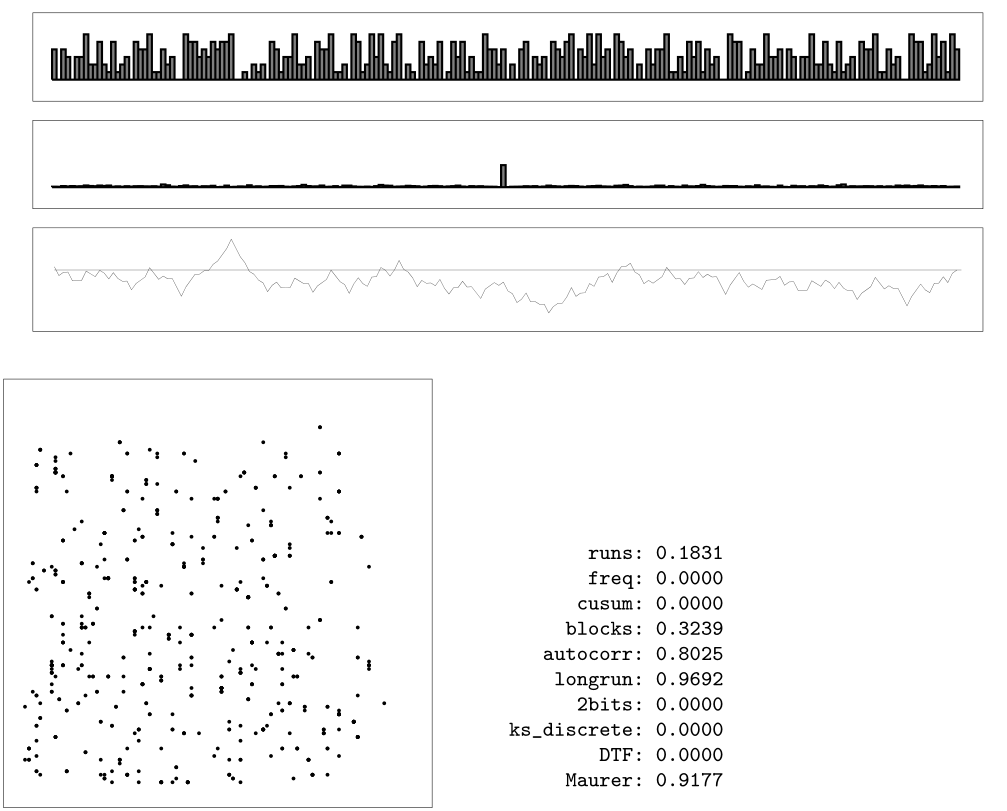}
\end{figure}

\newpage

\EDS{Example}{\rahmen{f(x,y) = (x+3y+3) \mod 4\\[2ex] a=3, b=2} }

\bgroup
\Century[950]
\begin{verbatim}
0
201
0210212
201201001201223
0210212210210030212210212232230
201201001201223201001201003003201201223201001201223223023223201
0210212210210030212210212232230210210030212210210030032030030210212 ...
2012010012012232010012010030032012012232010012012232230232232012010 ...
0210212210210030212210212232230210210030212210210030032030030210212 ...
2012010012012232010012010030032012012232010012012232230232232012010 ...
0210212210210030212210212232230210210030212210210030032030030210212 ...
2012010012012232010012010030032012012232010012012232230232232012010 ...
\end{verbatim}
\egroup
\begin{figure}[h]
	\includegraphics[width=1\linewidth]{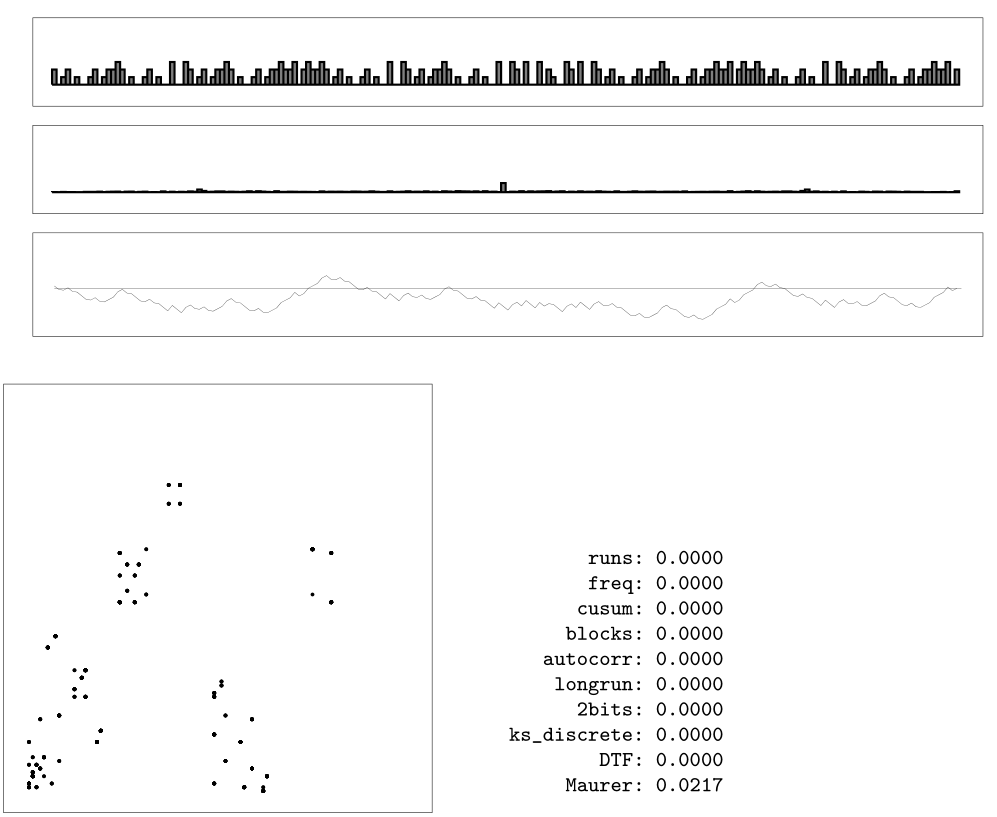}
\end{figure}

\newpage

\EDS{Example}{\rahmen{f(x,y) = (3x+5y+2) \mod 7\\[2ex] a=3, b=4}  }

\bgroup
\Century[950]
\begin{verbatim}
3
533
1543533
212524131543533
0261125562045163212524131543533
405236413112550556221014451106030261125562045163212524131543533
3400656223665451632131125505306505562242615001446445113150466033405 ...
5334002046355622422366163524451106030261632131125505306543404635306 ...
1543533400205210142603150556224204324223661641060315620464451131504 ...
2125241315435334002052106562615001443236603321253065055622420432101 ...
0261125562045163212524131543533400205210656261504635562236412530200 ...
4052364131125505562210144511060302611255620451632125241315435334002 ...
\end{verbatim}
\egroup
\begin{figure}[h]
	\includegraphics[width=1\linewidth]{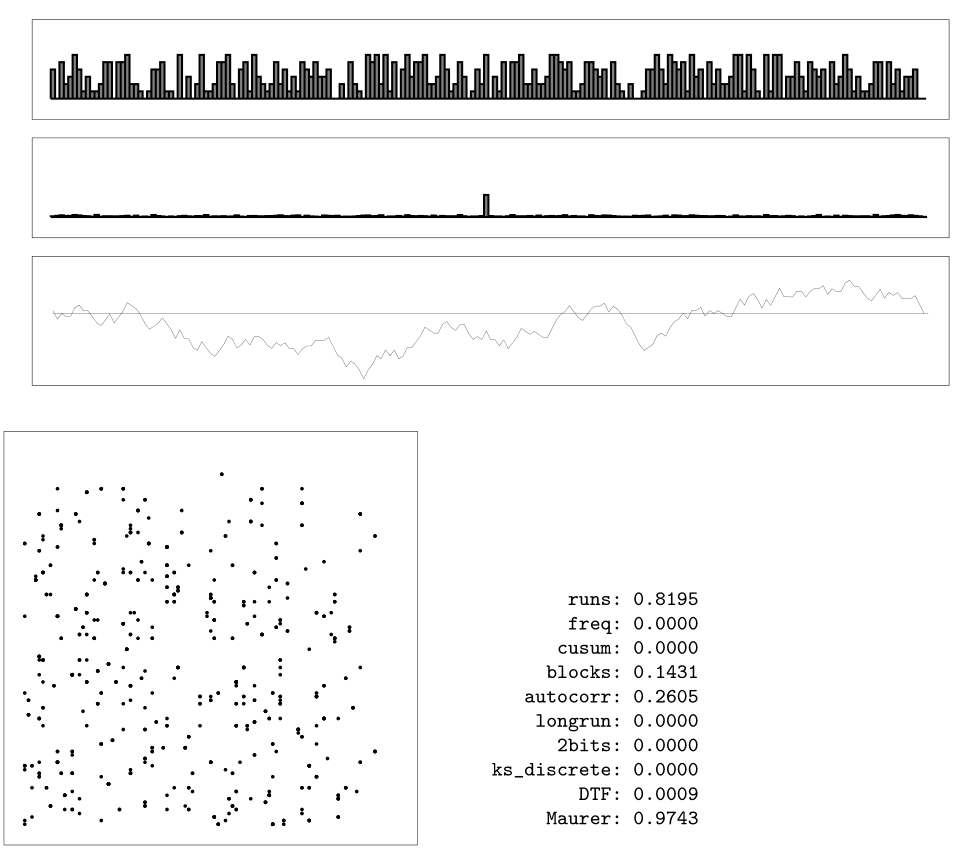}
\end{figure}

\newpage

\EDS{Example}{\rahmen{f(x,y) = (7x+4y) \mod 9\\[2ex] a=2, b=3}   }

\bgroup
\Century[950]
\begin{verbatim}
8
185
0138452
504113880435728
7580745121138878207443550732185
676548201724353162012113887837081250172484435515801773220138452
2677168564681250418732344355237146525041620121138878370853474058616 ...
4226775781462845564476286162758074513837732283148443551572834781547 ...
3402422677576507086154765218043515564484271652188611465267654820172 ...
8314108234024226775765071685801740588611056427168572013820744355310 ...
1853715451705812831410823402422677576507168580178146284548204187241 ...
0138452347810564353187402548616218537154517058128314108234024226775 ...
\end{verbatim}
\egroup
\begin{figure}[h]
	\includegraphics[width=1\linewidth]{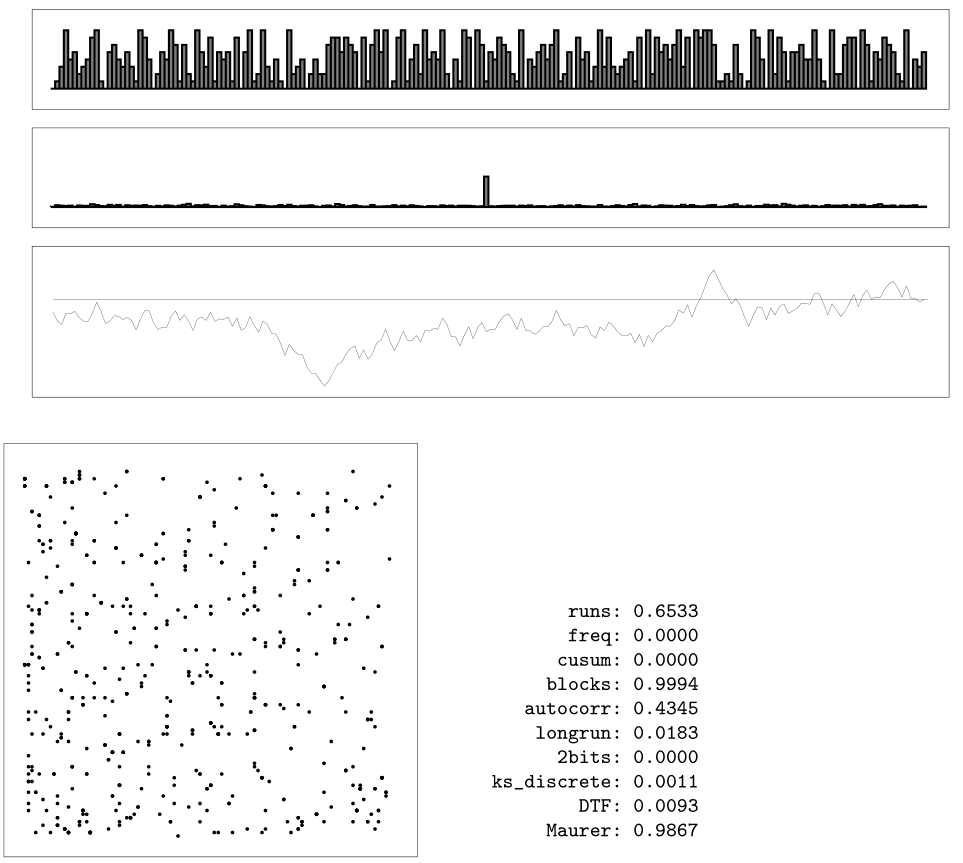}
\end{figure}

\newpage

\EDS{Example}{\rahmen{f(x,y) = (7x+4y+5) \mod 9\\[2ex]\hfill a=2, b=5}    }

\bgroup
\Century[950]
\begin{verbatim}
3
431
8403315
685460832331556
7678052436201813724323315565164
276637587075328403561210018821630782840372432331556516458106048
0227668653072508870067257372685460831516612251305001883862510653806 ...
1042022766867846457380678235401838870050263782355733078276780524362 ...
5130345210420227668678463758543604855733487026375862431524600188134 ...
3581638083648532513034521042022766867846375854365307250805240356203 ...
4315082106534870181356042805737235816380836485325130345210420227668 ...
8403315540186251302645736428870001882163151620345268707557330782431 ...
\end{verbatim}
\egroup
\begin{figure}[h]
	\includegraphics[width=1\linewidth]{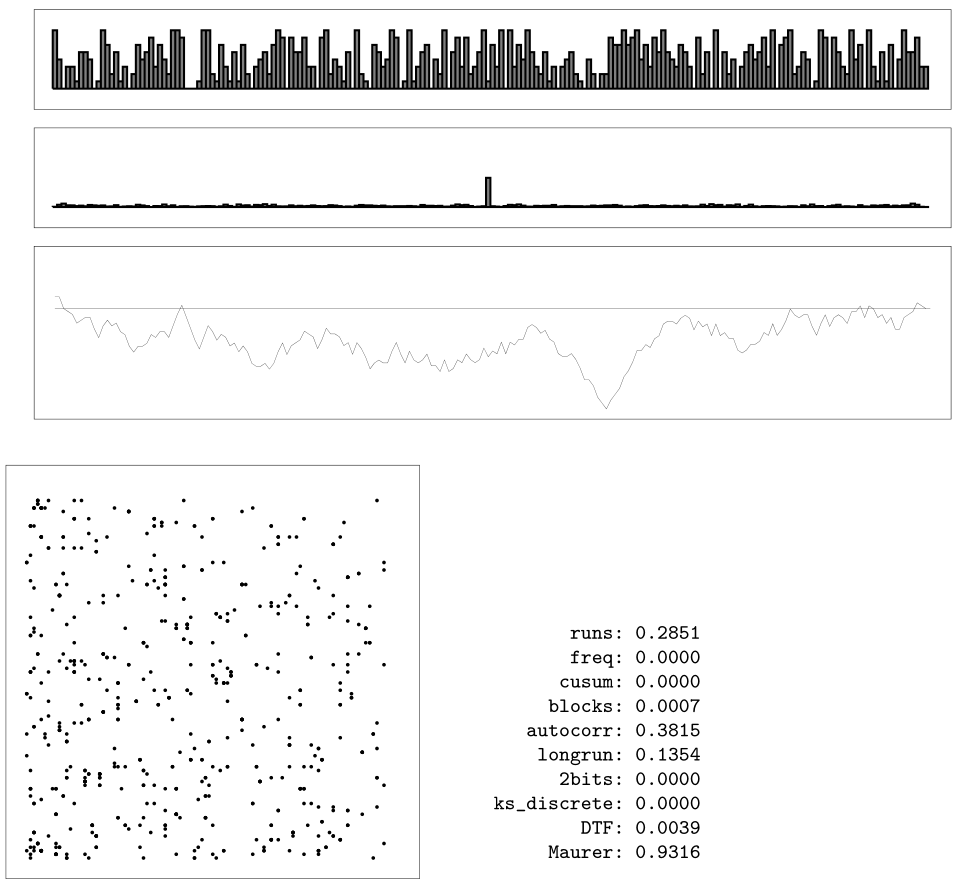}
\end{figure}

\newpage

\EDS{Example}{\rahmen[0]{f(x,y) = (x^3+2xy^2+x^2y+2y^3+5x^2+2xy+7y^2+6x+6y+7) \mod 9\\[2ex] a=1, b=8}    }

\bgroup
\Century[950]
\begin{verbatim}
8
883
8838236
883823686283766
8838236862837668860268231776766
883823686283766886026823177676683886700256686283713707764776766
8838236862837668860268231776766838867002566862837137077647767668236 ...
8838236862837668860268231776766838867002566862837137077647767668236 ...
8838236862837668860268231776766838867002566862837137077647767668236 ...
8838236862837668860268231776766838867002566862837137077647767668236 ...
8838236862837668860268231776766838867002566862837137077647767668236 ...
8838236862837668860268231776766838867002566862837137077647767668236 ...
\end{verbatim}
\egroup
\begin{figure}[h]
	\includegraphics[width=1\linewidth]{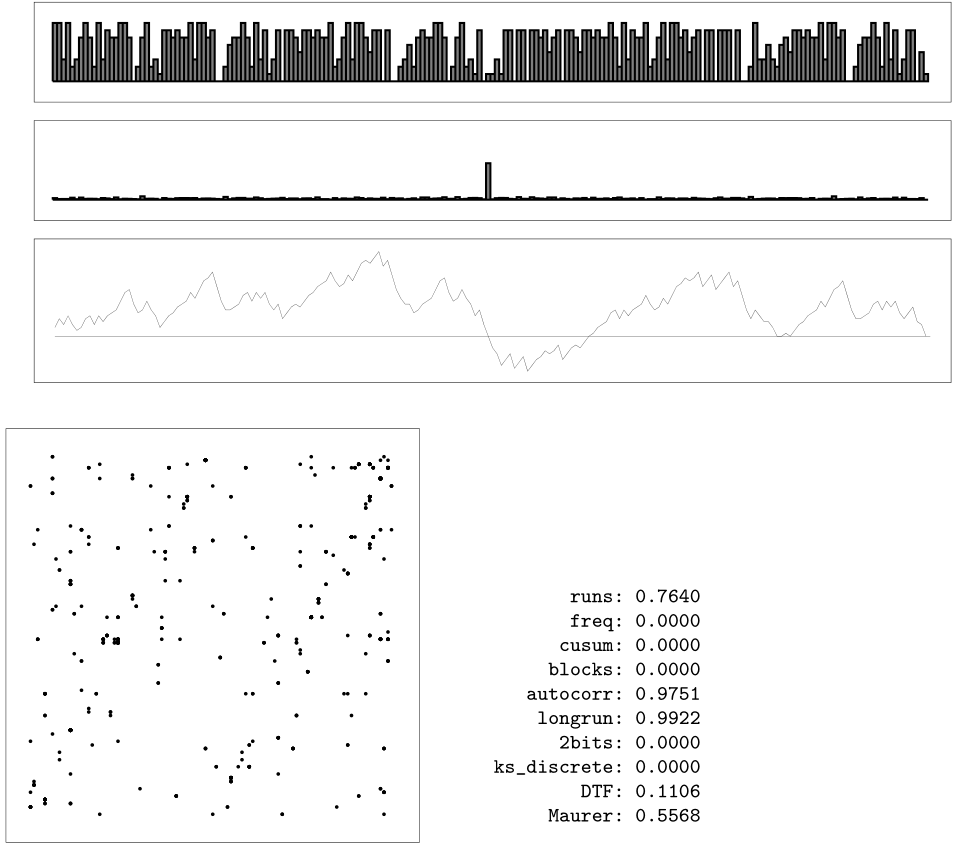}
\end{figure}

\newpage

\EDS{Example}{\rahmen[5]
	{f(x,y) = A^x_y, \quad \text{ where }
		A = \begin{pmatrix}
		1 & 4 & 2 & 5 & 3\\4 & 1 & 3 & 2 & 5\\5 & 2 & 4 & 3 & 1\\
			3 & 5 & 1 & 4 & 2\\2 & 3 & 5 & 1 & 4	
		\end{pmatrix}\\[3em]
		 a=1, b=4}    }

\bgroup
\Century[950]
\begin{verbatim}
5
351
2315215
423351353241351
5452334315212315532224312315215
351425323343341351353241423351354553221212241351423351353241351
2315215452555322334334134334312315212315532224315452334315212315142 ...
4233513532413514253255454553221233433413433431233413433413514233513 ...
5452334315212315532224312315215452555322554514251425455322124142334 ...
3514253233433413513532414233513545532212122413514233513532413514253 ...
2315215452555322334334134334312315212315532224315452334315212315142 ...
4233513532413514253255454553221233433413433431233413433413514233513 ...
\end{verbatim}
\egroup
\begin{figure}[h]
	\includegraphics[width=1\linewidth]{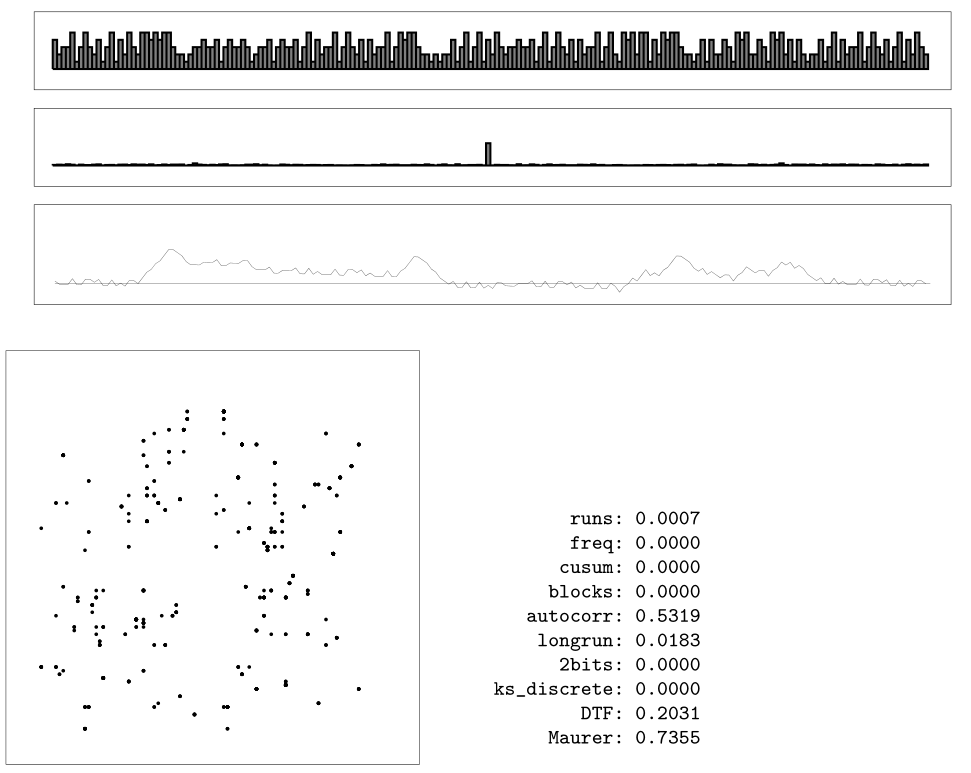}
\end{figure}

\newpage

\EDS{Example}{\rahmen[0]{f(x,y) = \Sistema {(3x+4y+1) \mod 9 & \text{ if } (x^2+y^3)\equiv 1\mod 8 \\ (7x+7y+4)\mod 9 & \text{ otherwise }}\\[3ex] a=3, b=4}     }

\bgroup
\Century[950]
\begin{verbatim}
8
687
7638172
271643682167321
0247611624837638227116572302712
100214071681011662147803271643682252476101163577323310024761321
5140400271245087611638214051011676627124070860730247611624837638225 ...
6511245054504002476132144570681716810116436822712450351140510116571 ...
7635110132144570355445705450400214071681530271246445774016382167611 ...
2716436511014051530271246445774073652554644577403554457054504002712 ...
0247611624837635110140512450351115631002476132142624644577372450872 ...
1002140716810116621478032716436511014051245035113214457073651101011 ...
\end{verbatim}
\egroup
\begin{figure}[h]
	\includegraphics[width=1\linewidth]{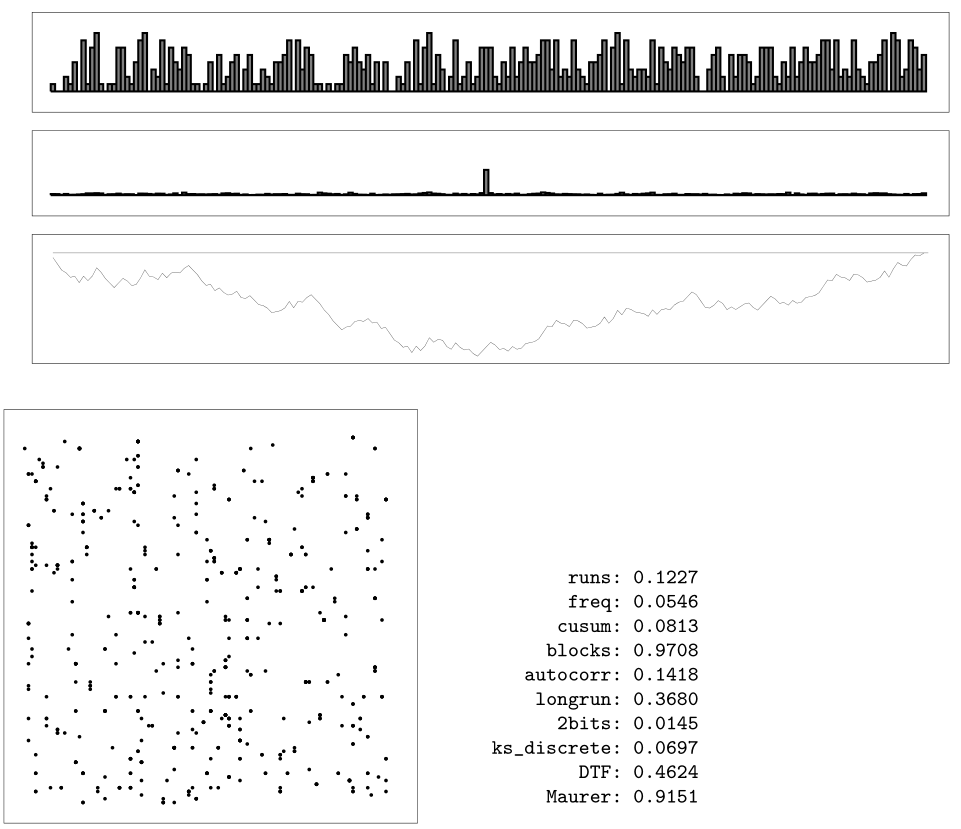}
\end{figure}

\newpage

\EDS{Example}{\rahmen{f(x,y) = (\on{altsum}(31x+35y+47)\mod 9 \\[2ex] a=18, b=11}
\sm\NI where $\on{altsum}(n)$ is the alternating sum of decimal digits of $n$.\m    }

\bgroup
\Century[950]
\begin{verbatim}
0
203
4210738
848211804763287
5864082211316870040716831248172
857816844088223211315381267817006004404721267843815274080147428
5865775801267864344088382232531211315381051328015226575801470060462 ...
8578163577571578706152265758168413343440883843282232531275138152113 ...
5865775801268375775715772105775817004641056232263577157801267864815 ...
8578163577571578706152267843071577571577210577574211802577571578014 ...
5865775801268375775715772105775817004641056232265758641360472105775 ...
8578163577571578706152267843071577571577210577574211802577571578014 ...
\end{verbatim}
\egroup
\begin{figure}[h]
	\includegraphics[width=1\linewidth]{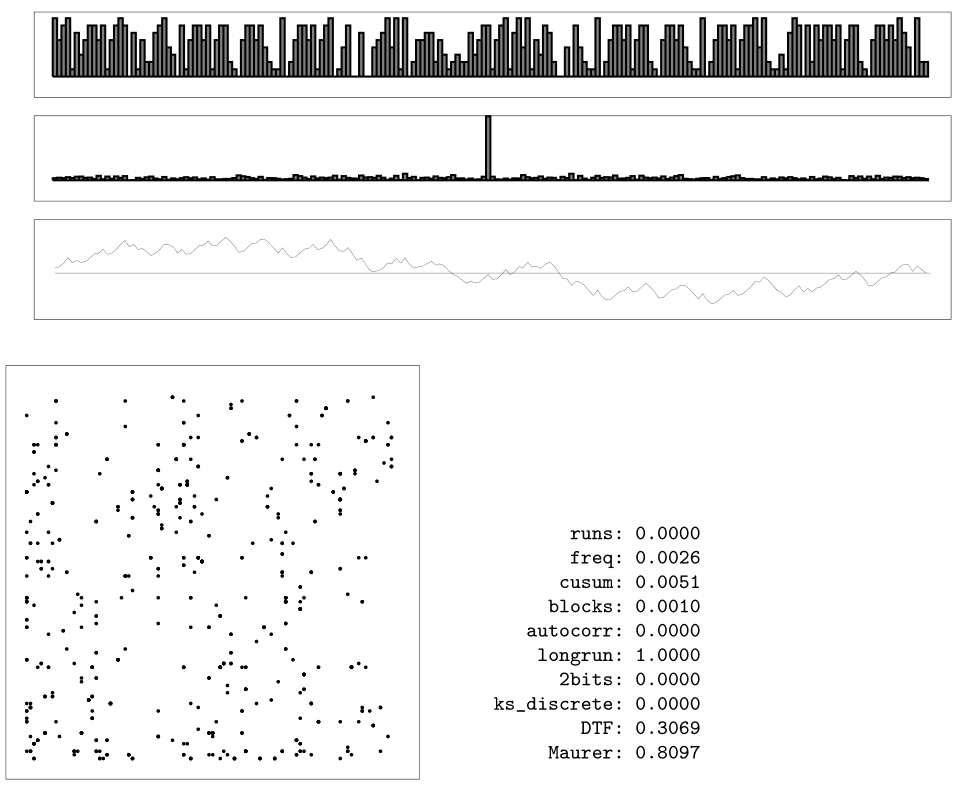}
\end{figure}

\newpage

\EDS{Example}{\rahmen{f(x,y) = \Mkr{\Mke{\Frac{x^2}{y}}+
			\Mke{\Frac{y^2}{x}}}\mod 7 +1\\[3ex] a=3, b=4}     }

\bgroup
\Century[950]
\begin{verbatim}
1
313
7331331
474373313373313
1417141347437331337347437331331
313431171134313314171413474373313373474314171413474373313373313
7331331413313117113133141331337331343117113431331417141347437331337 ...
4743733133733134313373313331311711313331337331343133733133734743733 ...
1417141347437331337347437331331413313373474373313373733133313117113 ...
3134311711343133141714134743733133734743141714134743733133733134313 ...
7331331413313117113133141331337331343117113431331417141347437331337 ...
4743733133733134313373313331311711313331337331343133733133734743733 ...
\end{verbatim}
\egroup
\begin{figure}[h]
	\includegraphics[width=1\linewidth]{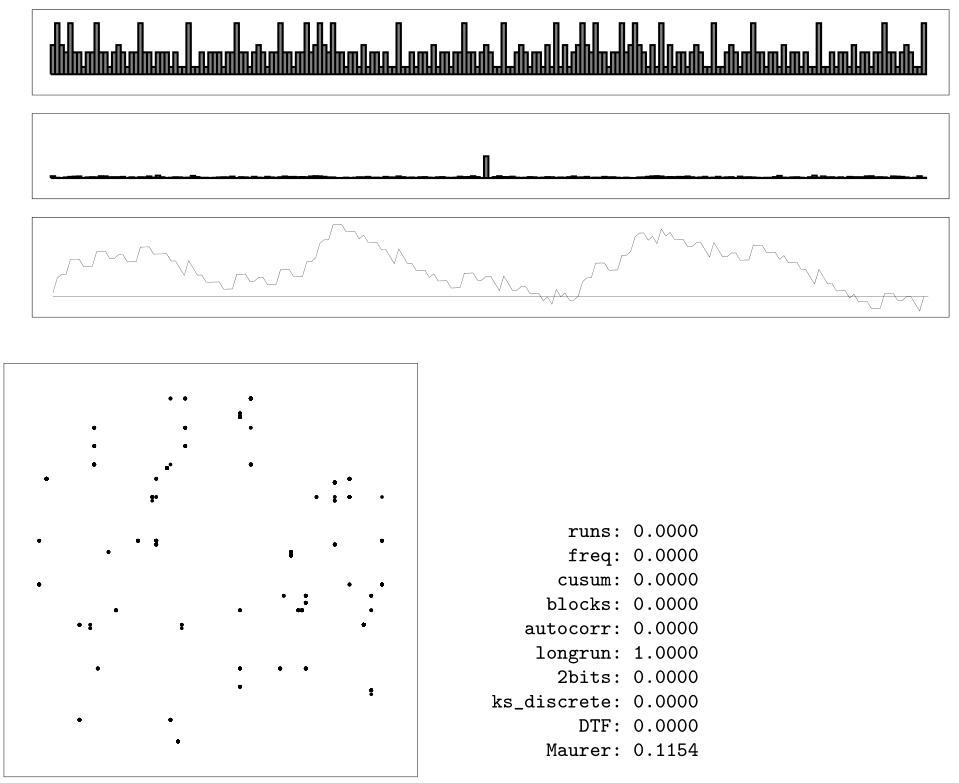}
\end{figure}

\newpage

\EDS{Example}{\rahmen{f(x,y) = \Mkr{\Mke{\Frac{x^2}{y+1}}+3}\mod 10\\[3ex] a=3, b=4}      }

\bgroup
\Century[950]
\begin{verbatim}
4
446
4464560
446456045566903
4464560455669034557566866940334
446456045566903455756686694033445575671566867826866994903353446
4464560455669034557566866940334455756715668678268669949033534464557 ...
4464560455669034557566866940334455756715668678268669949033534464557 ...
4464560455669034557566866940334455756715668678268669949033534464557 ...
4464560455669034557566866940334455756715668678268669949033534464557 ...
4464560455669034557566866940334455756715668678268669949033534464557 ...
4464560455669034557566866940334455756715668678268669949033534464557 ...
\end{verbatim}
\egroup
\begin{figure}[h]
	\includegraphics[width=1\linewidth]{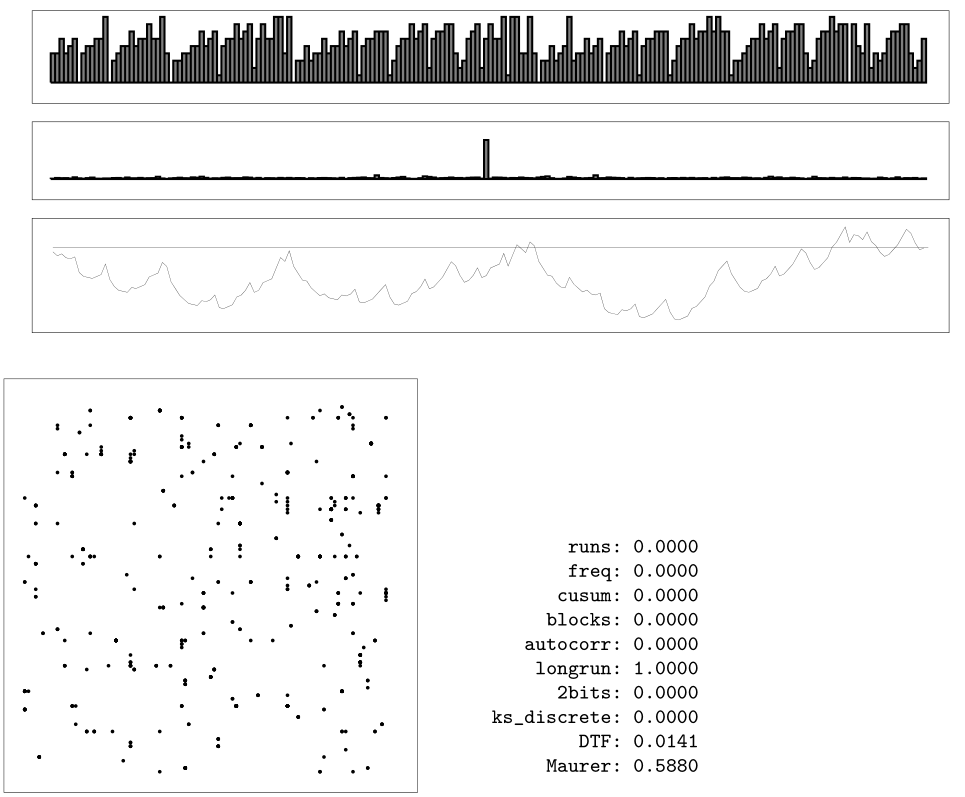}
\end{figure}

\newpage

\EDS{Example}{\rahmen{f(x,y) = (\on{gcd}(3x+4y+1,xy+y^2+4))\mod 5\\[2ex] a=3, b=4}       }

\bgroup
\Century[950]
\begin{verbatim}
2
221
2232114
223243221121441
2232432214032232112122114414114
223243221403223211441033223243221121221122321121441411441121441
2232432214032232114410332232432211214414114033232232432214032232112 ...
2232432214032232114410332232432211214414114033232232432214032232112 ...
2232432214032232114410332232432211214414114033232232432214032232112 ...
2232432214032232114410332232432211214414114033232232432214032232112 ...
2232432214032232114410332232432211214414114033232232432214032232112 ...
2232432214032232114410332232432211214414114033232232432214032232112 ...
\end{verbatim}
\egroup
\begin{figure}[h]
	\includegraphics[width=1\linewidth]{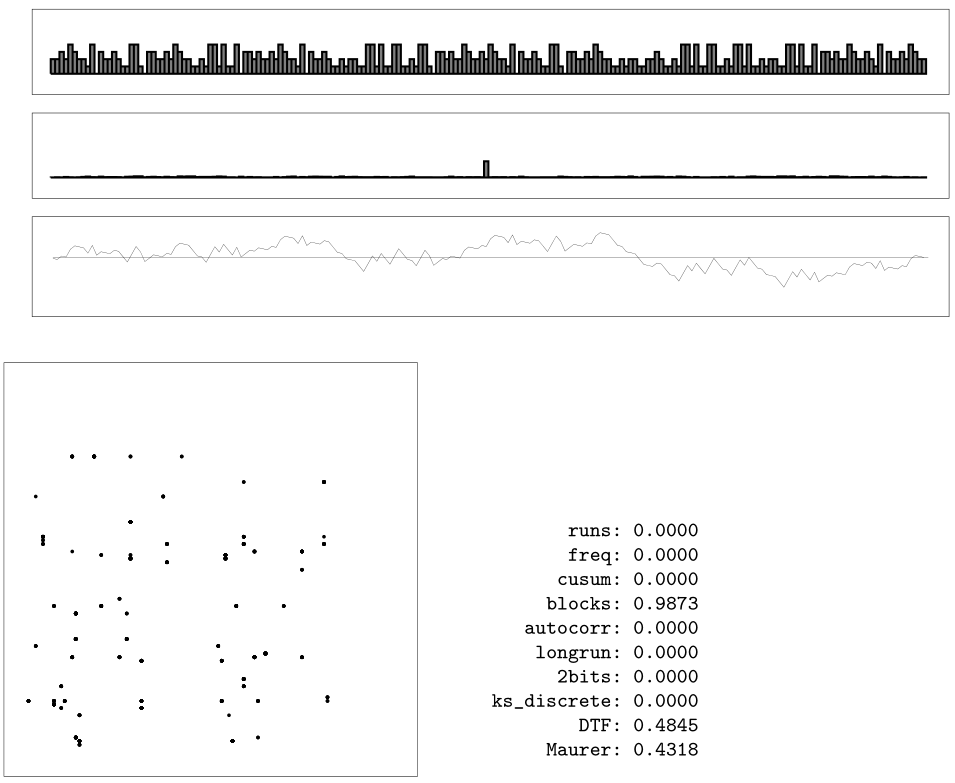}
\end{figure}

\newpage

\EDS{Example}{\rahmen{f(x,y) = |x-y+1|\\[2ex] a=2, b=7}    }

\bgroup
\Century[950]
\begin{verbatim}
4
142
2124324
122102142322142
2102122120122124320322122124324
122120122102122102300102122102142322302322122102122102142322142
2102122102300102122120122102122120120340100120122102122120122124320 ...
1221201221021221201203401001201221021221023001021221201221021221023 ...
2102122102300102122120122102122102300102302304500120100102300102122 ...
1221201221021221201203401001201221021221023001021221201221021221201 ...
2102122102300102122120122102122102300102302304500120100102300102122 ...
1221201221021221201203401001201221021221023001021221201221021221201 ...
\end{verbatim}
\egroup
\begin{figure}[h]
	\includegraphics[width=1\linewidth]{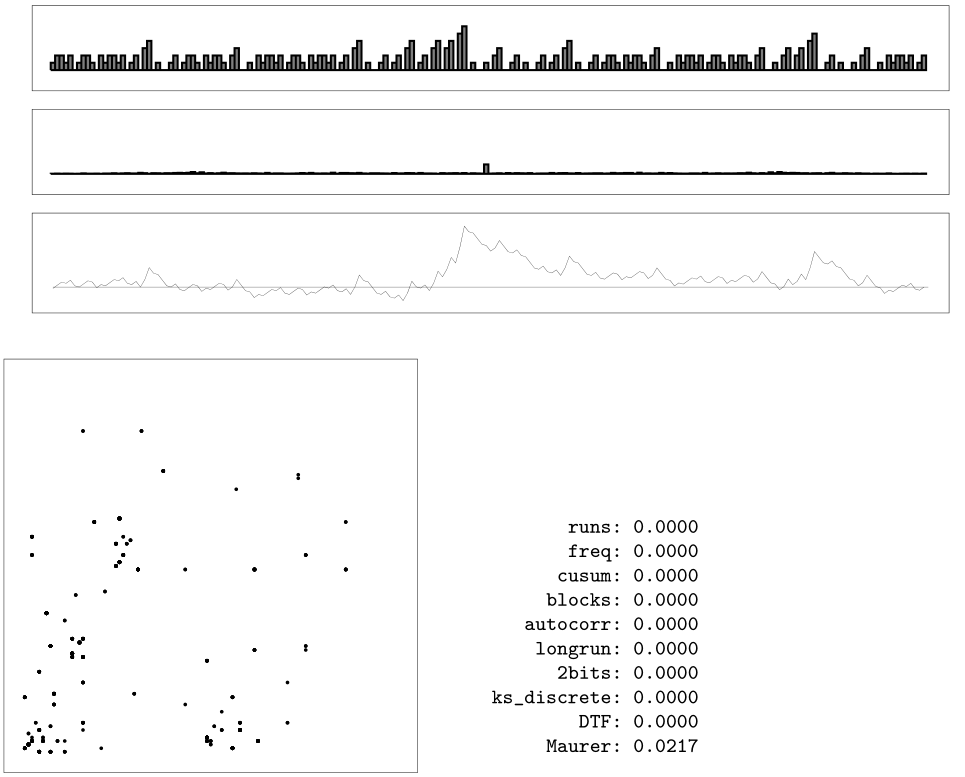}
\end{figure}

\Bibliografia[References]{

\fonte[0] {25668} {L. Bassham a.o.} {A statistical test suite
for random and pseudorandom number generators for cryptographic applications}.
NIST 2010.

\fonte[1] {9797} {R. Centrella}  {Numeri casuali - teoria e generatori dicotomici}. \\
Tesi Univ. Roma 1997.

\fonte[1] {26401} {J. Coron/D. Naccache} {An accurate evaluation of Maurer's universal test}.
Springer LN CS 1556 (2002), 57-71.

\fonte[1] {26390} {A. Do\u{g}anaksoy/C. Tezcan} {An alternative approach to Maurer's universal
statistical test}. \\
3rd Inf. Sec. Crypt. Conference Ankara 2008.

\fonte[1] {2031} {M. Fisz} {Wahrscheinlichkeitsrechnung und mathematische Statistik}.\\
Deutscher Vlg. Wiss. 1989.

\fonte[1] {25885} {L. Guibas/A. Odlyzko} {Long repetitive patterns in random sequences}.
Zt. Wtheorie verw. Geb. 53 (1980), 241-262.

\fonte[1] {12} {D. Knuth} {The art of computer programming}.\\ Volume 1. Addison-Wesley.

\fonte[1] {24739} {H. Kreindl} {BUS-Theorie}. Internet ca. 2012.

\fonte[1] {19061} {L. Kuipers/H. Niederreiter} {Uniform distribution of sequences}. \\
Dover 2006.

\fonte[1] {25830} {U. Maurer} {A universal statistical test for random bit generators}. \\
J. Cryptology 5/2 (1992), 89-105.
}

\end{document}